\documentclass[preprint,12pt]{elsarticle}

\usepackage{amsthm}
\usepackage{amssymb}
\usepackage{epsfig}
\usepackage{graphicx}
\usepackage{epstopdf}
\usepackage{amsthm}
\usepackage{srcltx}
\usepackage[dvipsnames,usenames]{color}
\usepackage{paralist}
\usepackage{amsmath}
\usepackage{amsfonts}

\usepackage{multicol}
\usepackage{bm}
\usepackage{float}
\usepackage{graphicx}
\usepackage{caption}
\usepackage{subcaption}
\usepackage{abraces}
\usepackage{amscd}
\usepackage{amstext}
\usepackage{color}
\usepackage{textcomp}
\usepackage{stmaryrd}
\usepackage{pstricks}

\usepackage{pgfplots}

\DeclareMathAlphabet{\pazocal}{OMS}{zplm}{m}{n}

\newcounter{contador}
\newcounter{teoA}

\newtheorem{teoa}[teoA]{Theorem}

\newtheorem{theorem}[contador]{Theorem}

\newtheorem{definition}[contador]{Definition}

\newcounter{ex}

\begin{document}

\begin{frontmatter}

\title{Minor loops of the Dahl and LuGre models}

\author[a1]{Fay\c{c}al Ikhouane\corref{mycorrespondingauthor}} \cortext[mycorrespondingauthor]{Corresponding author}
\author[a2]{V\'{\i}ctor Ma\~{n}osa}
\author[a2]{Gisela Pujol}

\address{Departament de Matem\`{a}tiques, Universitat Polit\`{e}cnica de Catalunya}

\address[a1]{EEBE, Av. Eduard Maristany 16, 08019 Barcelona, Spain}
\address[a2]{ESEIAAT, Colom 1 and 11, 08222 Terrassa, Spain}

\begin{abstract}
Hysteresis is a special type of behavior encountered in physical systems: in a hysteretic system, when the input  is periodic and varies slowly, the steady-state part of the  output-versus-input  graph becomes a loop called \emph{hysteresis loop}. In the presence of  perturbed inputs, this hysteresis loop presents small lobes called minor loops that are located inside a larger curve called major loop. The study of minor loops is being increasingly popular since it leads to a quantification of the loss of energy. The aim of the present paper is to give an explicit analytic expression of the  minor loops of the LuGre and the Dahl models of dynamic dry friction. 
\end{abstract}

\begin{keyword}
Hysteresis; Minor loops; LuGre and Dahl models.

\noindent {\sl  MSC 2010:} 34C55; 93A30; 93A99; 46T99; 

\noindent {\sl  PACS:} 77.80.Dj; 75.60.-d.
\end{keyword}

\end{frontmatter}

\section{Introduction}

Hysteresis is a nonlinear phenomenon observed in some physical systems under low-frequency excitations. It appears in many areas as biology, electronics, ferroelasticity, magnetism, mechanics or optics \cite{BM2006,BS1996,JMS2017,MNZ93,V94}. This phenomenon is currently  classified into two categories: \emph{rate independent} (RI) and \emph{rate dependent} (RD) hysteresis. For RI hysteresis, the output-versus-input graph of the hysteresis system does not change with the frequency of the input signal. This is the case for example of the Bouc-Wen or the Preisach models, see \cite{IR2007} and \cite{M03} respectively. For RD hysteresis, the output-versus-input graph of the hysteresis system may change with the frequency, but it converges in some sense to a fixed loop called the hysteresis loop when the frequency goes to zero. This is the case for example of the LuGre model and the semilinear Duhem model, see  \cite{Naser et al.(2015)} and \cite{Ikhouane2017,OB05} respectively.    Research in the field of hysteresis has focused mainly on the study of rate-independent hysteresis, and it is only in the last 15 years that the importance of rate-dependent phenomena has been acknowledged, and it constitutes a challenge by itself. 

\begin{figure}[H]
                \centering
                 \includegraphics[scale=0.15]{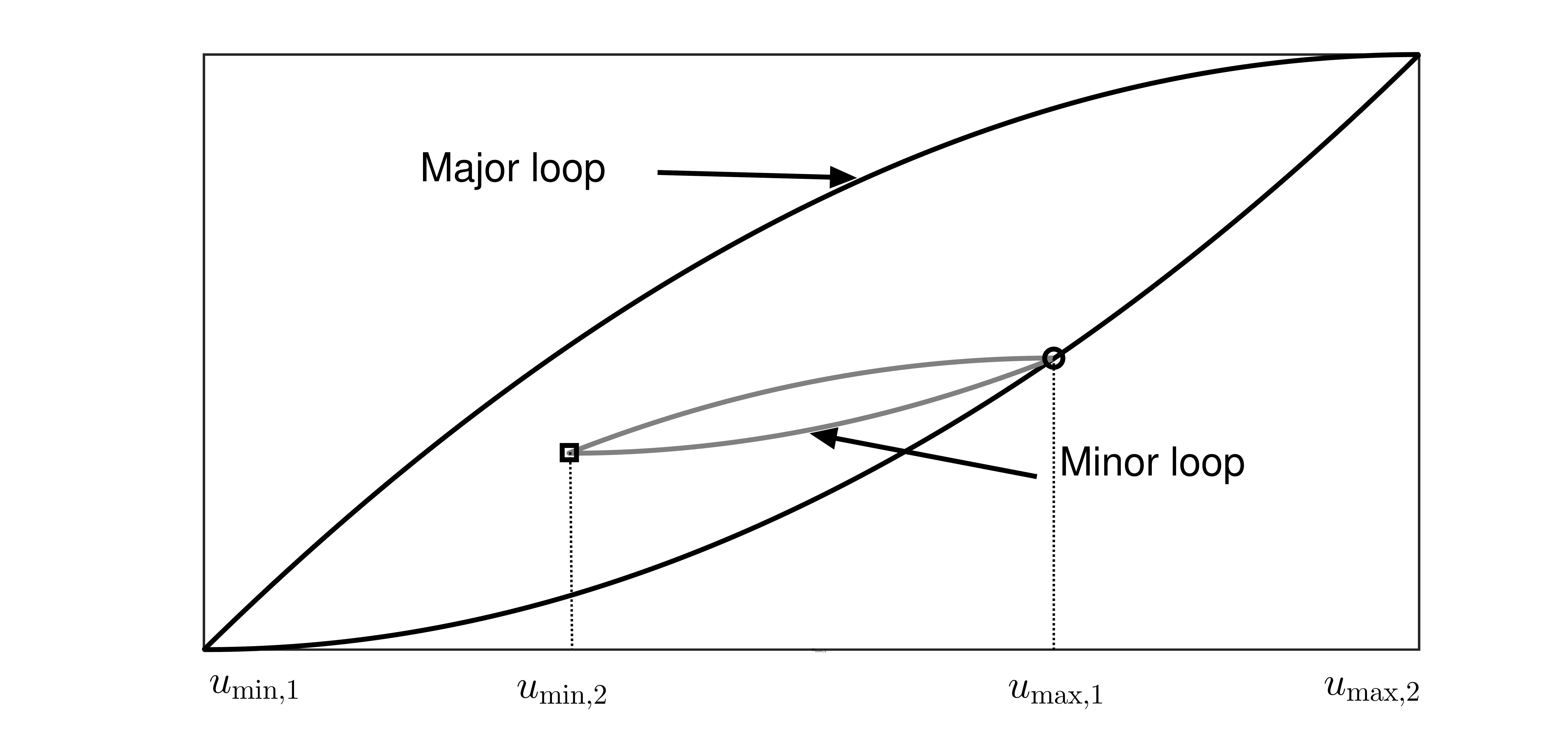}
                \caption{Hysteresis loop of a magnetic system with the input of Figure  \ref{f:figugen}.  Black: major loop. Grey: minor loop.}
                \label{plot66yX5y}
\end{figure}

\begin{figure}[H]
                \centering
                \includegraphics[scale=0.15]{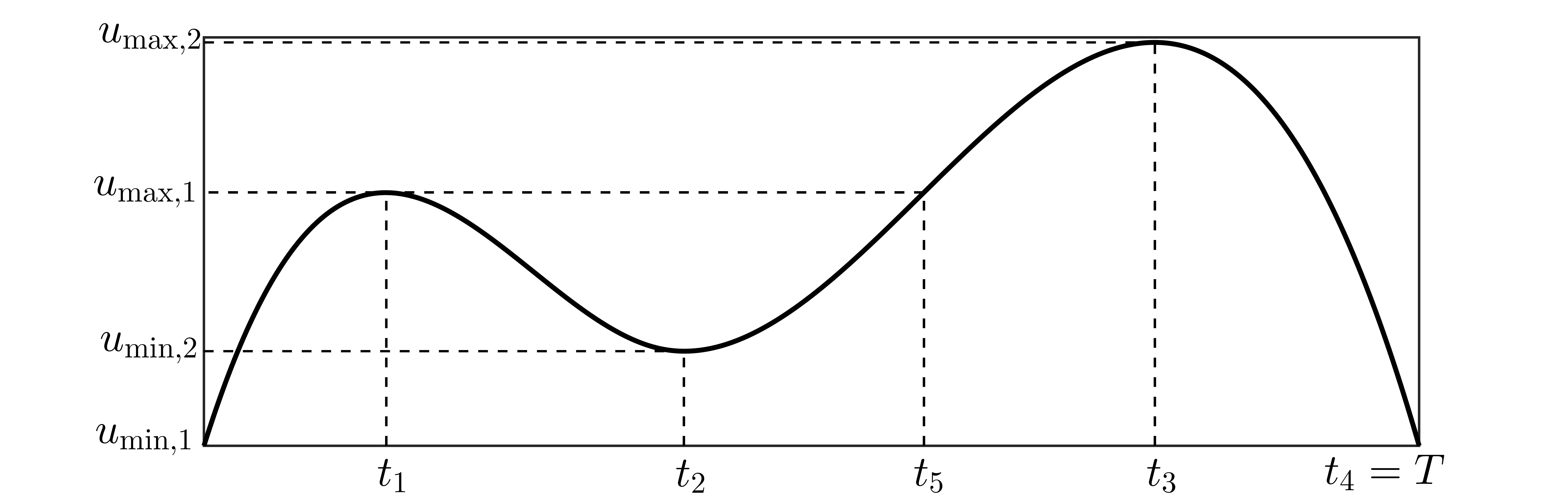}
                \caption{A bimodal $T$--periodic input $u(t)$ versus $t$.}
                \label{f:figugen}
\end{figure}

The recent years have witnessed a growing interest in a phenomenon that appears in hysteretic systems under perturbed  periodic signals: the hysteresis loop shows to be composed of a big cycle called major loop, and one or several small lobes called minor loops located inside the major loop. Figure \ref{plot66yX5y} shows the hysteresis loop of a magnetic system when the input is the one of Figure \ref{f:figugen}, see  \cite{HMFA14} for instance.

 This interest in the study of minors loops is due in part  to the fact that minor loops are related to a loss of energy, see  \cite{ZRdKAF2017} for instance.

From a formal point of view, minor loops have been studied mainly in relation with the Preisach model \cite[p.19]{M03}. Apart from \cite[Sections 10, 11.9]{Ikhouane2017} we are not aware of any mathematical analysis of minor loops of hysteresis systems given by differential equations. 

The aim of the present paper is to fill this void by providing an explicit analytic description of the minor loops of the Dahl and the LuGre models.

The Dahl model is an idealization of dynamic dry friction proposed by Dahl in 1976 \cite{D1976}. This model relates an input displacement $u$ to an output force $y$ as
\begin{align*}
&y(t)=F_{c}\,w(t), \\
&\dot{w}(t)=\rho\big(\dot{u}(t)-|\dot{u}(t)|w(t)\big), \\
&-1 \leq w (0) \leq 1,  
\end{align*}
where $w$ is an internal state and $\rho>0$, $F_c>0$ are constants. A good introductory text on the relationship between the Dahl model and the Coulomb model of dry friction may be found in \cite{GBI2016}.

 The LuGre model is a generalization of the Dahl model introduced in 1995 to include the \emph{Stribeck effect},   that is the decrease of friction at low velocities \cite{wc95}. The LuGre model is  given by \cite{JC08}:
\begin{equation}\label{equation1}
\begin{array}{lcl}
\dot{x}\left(t\right) & = & -\sigma_{0}\dfrac{|\dot{u}(t)|}{g\big(\dot{u}(t)\big)}x\left(t\right)+\dot{u}(t),  \\
x(0) & = & x_{0},  \\
F\left(t\right)  & = & \sigma_{0}x(t)+\sigma_{1}\dot{x}(t)+f\big(\dot{u}(t)\big), 
\end{array}
\end{equation}

\noindent where $t\geq0$ denotes time; the parameters $\sigma_{0}>0$ and $\sigma_{1}>0$ are respectively the stiffness and the microscopic damping friction coefficients; the function $g$ is continuous with $g\left(\vartheta\right)>0$  for all $\vartheta\in\mathbb{R}$ and it represents the macrodamping friction ; $x(t)\in\mathbb{R}$ is the average deflection of the bristles; $x_{0}\in\mathbb{R}$ is the initial state; $u(t)$ is the relative displacement and is the input of the system; $F(t)$ is the friction force and is the output of the system; and $f$ is continuous and such that $f(0)=0$.  
When the function $g$ is constant, $\sigma_1=0$ and $f$ is the zero function, the system (\ref{equation1}) reduces to  the Dahl model. Both the LuGre and the Dahl models have been used in various applications, see for instance \cite{AIRC2012,FHH2008,GP2017}.

The main contribution of this paper is Theorem \ref{t:main} which is stated in Section \ref{ss:statement}. This theorem provides  the analytic description of  the minor loop of the LuGre and  Dahl models when the input is bimodal like in Figure \ref{f:figugen}.

The paper is organized as follows:  Section \ref{notation} provides  the mathematical notation used in the text.
In Section \ref{s:main} we present and prove the main result which is the analytical description of the  minor loop of the LuGre and Dahl models.   Section \ref{s:numsim} has a pedagogical interest:  using  numerical simulations  we present  examples that illustrate the constructive process which leads to the hysteresis and minor loops. Some of these examples are 
 aimed for the reader who may not be familiar with the technicalities that underline the methodology used here. The conclusions are provided in Section \ref{Conclusions and comments}.

\section{Mathematical notation}\label{notation}

 We say that a subset of $\mathbb{R}$ is measurable when it is Lebesgue measurable. Consider a function $g: I \subseteq \mathbb{R} \rightarrow \mathbb{R}$ where  $I $ is an interval.
We say that $g$ is measurable if $g^{-1}(B)$ is a measurable set for any set $B$ in the Borel algebra of $\mathbb{R}$ or, equivalently, if $\{x\in I:\, g(x)>a\}$ is a measurable set for all $a\in\mathbb{R}$, \cite{Rudin, YV1966}. For a measurable function $g:I\rightarrow \mathbb{R}$, $\|g\|$ denotes the essential supremum of the function $|g|$ where $|\cdot|$ is the absolute value. 

We recall that $\mathcal{C}^0(\mathbb{R}^+,\mathbb{R})$ denotes the space of continuous functions defined from $\mathbb{R}+$ to $\mathbb{R}$ endowed with the norm $\| \cdot\|$. Also
$W^{1;\infty}(\mathbb{R}^+,\mathbb{R})$ denotes the Sobolev space of  absolutely continuous functions $u:\mathbb{R}^+ \rightarrow \mathbb{R}$. For this class of functions, the derivative $\dot{u}$ is measurable, and we have  $\|u\| <\infty$, $\|\dot{u}\|<\infty$. Endowed with the norm $\|u\|_{1,\infty}=\max\left(\|u\|, \|\dot{u}\|\right)$, $W^{1;\infty}(\mathbb{R}^+,\mathbb{R})$ is a Banach space \cite[pp.\,280--281]{Leoni09}. Finally, $L^\infty(I,\mathbb{R})$ denotes the Banach space of measurable functions $u : I \rightarrow \mathbb{R}$ such that $\|u\| <\infty$, endowed with the norm $\| \cdot \|$. For  $T>0$ we define $\Omega_T$ as the set of all $T$--periodic functions~$u\in W^{1;\infty}(\mathbb{R}^+,\mathbb{R})$.

\section{Main result}\label{s:main}

\subsection{Statement of the main result}\label{ss:statement}
We consider the LuGre model  (\ref{equation1}) with an input $u \in W^{1;\infty}(\mathbb{R}^+,\mathbb{R})$.
In \cite{Naser et al.(2015)} it is proved that for all $x_0 \in \mathbb{R}$, the differential equation (\ref{equation1}) has a unique Carath\'eodory solution $x \in W^{1;\infty}(\mathbb{R}^+,\mathbb{R})$  and that $F \in L^\infty(\mathbb{R}^+, \mathbb{R})$.

To present the main result of this work which is the analytic characterization of  the minor loop we define formally the set of bimodal inputs needed to generate this minor loop.

\begin{definition}\label{a:ass1} Let $u_{\min,1},u_{\min,2},u_{\max,1}, u_{\max,2} \in \mathbb{R}$ be such that $u_{\min,1} \leq u_{\min,2}<u_{\max,1}\leq u_{\max,2}$ and at least one of the following holds: $u_{\min,1} \neq u_{\min,2}$ or $u_{\max,1} \neq u_{\max,2}$. Let $t_1,t_2,t_3,t_4 \in \mathbb{R}^+$ be such that $0<t_1<t_2<t_3<t_4$. We define  $\mathbb{M}_{u_{\min,1},u_{\min,2},u_{\max,1}, u_{\max,2},t_1,t_2,t_3,t_4}$ as the set of all functions $u \in \Omega_{t_4}$ such that $u$ is strictly increasing on the interval $[0,t_1]$, strictly decreasing on the interval $[t_1,t_2]$, strictly increasing on the interval $[t_2,t_3]$, strictly decreasing on the interval $[t_3,t_4]$; and $u(0)=u_{\min,1}$, $u(t_1)=u_{\max,1}$, $u(t_2)=u_{\min,2}$, $u(t_3)=u_{\max,2}$, $u(t_4)=u(0)$.
\end{definition}

\begin{teoa}\label{t:main}
 Let us consider the LuGre model given by Equations \eqref{equation1} with an  input $u \in \mathbb{M}_{u_{\min,1},u_{\min,2},u_{\max,1}, u_{\max,2},t_1,t_2,t_3,t_4}$. Then the following statements hold:
\begin{enumerate}
\item[\textnormal{(a)}] The hysteresis  loop that corresponds to the input $u$ is the set 
$$
G_u^\circ = \big \{ \big(\psi_u(t),y^\circ(t)\big)\in\mathbb{R}^2,\,t \in [0,\varrho_4] \big \},
$$
where $y^\circ$ is given by
$$
y^\circ(t)=\mathrm{e}^{-\frac{\sigma_0}{g(0)}(t-\varrho_i)}\left( y^\circ(\varrho_i)- g(0)\left[ \mathrm{e}^{\frac{\sigma_0}{g(0)}(t-\varrho_i)} -1\right]\right), \mbox{ for  } t \in [\varrho_i,\varrho_{i+1}]
$$
 and  $i\in\{0,1,2,3\}$, and where
\begin{align*}
&y^\circ(0) =g(0)\,\frac{\mathrm{e}^{-\frac{\sigma_0}{g(0)}\varrho_4}}{1-\mathrm{e}^{-\frac{\sigma_0}{g(0)}\varrho_4}}\left(  2\mathrm{e}^{\frac{\sigma_0}{g(0)}\varrho_1} -2\mathrm{e}^{\frac{\sigma_0}{g(0)}\varrho_2}+2\mathrm{e}^{\frac{\sigma_0}{g(0)}\varrho_3}-\mathrm{e}^{\frac{\sigma_0}{g(0)}\varrho_4}-1\right) \mbox{ and } \\
&y^\circ(\varrho_i)=\mathrm{e}^{-\frac{\sigma_0}{g(0)}(\varrho_i-\varrho_{i-1})}\left( y^\circ(\varrho_{i-1})+ g(0)\left[ \mathrm{e}^{\frac{\sigma_0}{g(0)}(\varrho_i-\varrho_{i-1})} -1\right]\right) \mbox{ for }  i\in\{1,2,3,4\};
\end{align*}
and $\psi_u$ is given by $$\psi_u(t)=\begin{cases}\begin{array}{ll}
t+u_{\min,1}& \text{for } t\in[0,\varrho_1],\\
-t+\varrho_1+u_{\max,1}& \text{for } t\in[\varrho_1,\varrho_2],\\
t-\varrho_2+u_{\min,2}& \text{for } t\in[\varrho_2,\varrho_3],\\
-t+\varrho_3+u_{\max,2}& \text{for } t\in[\varrho_3,\varrho_4],
\end{array}
\end{cases}$$
being $\varrho_0=0$, $\varrho_1=\;u_{\max,1}-u_{\min,1}>0,$
$\varrho_2=\;\varrho_1+u_{\max,1}-u_{\min,2}>\varrho_1,$
$\varrho_3=\;\varrho_2+u_{\max,2}-u_{\min,2}>\varrho_2,$ and
$\varrho_4=\;\varrho_3+u_{\max,2}-u_{\min,1}>\varrho_3$.
\item[\textnormal{(b)}] The minor loop that corresponds to the input $u$ is  
 the set $$\mathcal{N}_u=\big\{\big(\psi_u(t),y^\circ(t)\big),t \in [\varrho_1,\varrho_5]\big\},$$ where $\varrho_5=u_{\max,1}-u_{\min,2}+\varrho_2 \in (\varrho_2,\varrho_3]$.
\end{enumerate}

\end{teoa} 

\textbf{Comment.} Observe that the sets $G^{\circ}_{u}$ and $\mathcal{N}_u$ are the geometric loci of  parametrized curves. Theorem \ref{t:main}, thus, gives an explicit parametrization of these curves.

\subsection{Proof of Theorem \ref{t:main}} \label{ss:proof}

The proof of Theorem \ref{t:main} is done is three steps: 
\begin{enumerate}
\item[\textbf{Step 1:}] The hysteresis loop of the LuGre and the Dahl models are derived in Section \ref{ss:hystloop}.
\item[\textbf{Step 2:}] A normalized input is presented in Section \ref{ss:normalizedinput}.
\item[\textbf{Step 3:}] The determination of the equations of the minor loop is done in Section \ref{Analytical expression of the hysteresis and minor loops}.
\end{enumerate}

\subsubsection{Hysteresis loop of the LuGre and the Dahl models}\label{ss:hystloop}

To prove Theorem \ref{t:main} and therefore  to derive the explicit expression of the hysteresis loop of the LuGre and the Dahl models, 
we follow the methodology presented in \cite{I09,Naser et al.(2015)}. In this section we recall and adapt the main steps of this methodology. The reader unfamiliar with this theoretical framework is first referred to Example 1 in Section \ref{ss:ex1}.

Let $u \in \mathbb{M}_{u_{\min,1},u_{\min,2},u_{\max,1}, u_{\max,2},t_1,t_2,t_3,t_4}$ and take $T=t_4$. Also, take $\gamma  \in (0,\infty)$ and consider the linear time-scale change $s_\gamma:\mathbb{R} \rightarrow \mathbb{R}$ defined by $s_\gamma(t)=t/\gamma$ for all $t \in \mathbb{R}$. Then $u \circ s_\gamma$ is $\gamma T$--periodic. 

The system 
(\ref{equation1}) for which the input is $u \circ s_\gamma$ can be written as
\begin{eqnarray*} 
\dot{x}_\gamma\left(t\right) & = & -\sigma_{0}\frac{\big|\dot{\aoverbrace[L1R]{{u \circ s_\gamma}}}(t)\big|}{g\big(\dot{\aoverbrace[L1R]{u \circ s_\gamma}}(t)\big)}x_\gamma\left(t\right)+\dot{\aoverbrace[L1R]{{u \circ s_\gamma}}}\left(t\right), \text{ for almost all } t \in \mathbb{R}^+,\\
x_\gamma(0) & = & x_{0}, \\
F_\gamma\left(t\right)  & = & \sigma_{0}x_\gamma\left(t\right)+\sigma_{1}\dot{x}_\gamma(t)+f\big(\dot{\aoverbrace[L1R]{{u \circ s_\gamma}}}(t)\big), \text{ for almost all } t \in \mathbb{R}^+. 
\end{eqnarray*}
On the other hand, $$\dot{\aoverbrace[L1R]{{u \circ s_\gamma}}}(t)=\dot{u} \big( s_\gamma(t)\big) \cdot \dot{s}_\gamma(t)=\dot{u}\left(\frac{t}{\gamma}\right) \cdot \frac{1}{\gamma},$$
so that we get
\begin{eqnarray*} 
\dot{x}_\gamma\left(t\right) & = & -\sigma_{0}\frac{\left| \frac{1}{\gamma}\dot{u}\left(\frac{t}{\gamma}\right) \right|}{g\left(\frac{1}{\gamma}\dot{u}\left(\frac{t}{\gamma}\right)\right)}x_\gamma\left(t\right)+\frac{1}{\gamma}\dot{u}\left(\frac{t}{\gamma}\right), \text{ for almost all } t \in \mathbb{R}^+,\\
x_\gamma(0) & = & x_{0}, \\
F_\gamma\left(t\right)  & = & \sigma_{0}x_\gamma\left(t\right)+\sigma_{1}\dot{x}_\gamma(t)+f\left(\frac{1}{\gamma}\dot{u}\left(\frac{t}{\gamma}\right)\right), \text{ for almost all } t \in \mathbb{R}^+. 
\end{eqnarray*}
Now, defining $\nu=t / \gamma $ we rewrite these relations in terms of $\nu$ obtaining
\begin{equation} \label{ecuacion1explicativa} \begin{split}
\gamma \dot{x}_\gamma\left(\gamma \nu \right) & =  -\sigma_{0}\frac{| \dot{u}(\nu) |}{g\left(\frac{1}{\gamma}\dot{u}(\nu)\right)}x_\gamma(\gamma \nu)+\dot{u}(\nu), \text{ for almost all } \nu \in \mathbb{R}^+,\\
x_\gamma(0) & =  x_{0}, \\
F_\gamma(\gamma \nu)  & =  \sigma_{0}x_\gamma(\gamma \nu)+\sigma_{1}\dot{x}_\gamma(\gamma \nu)+f\left(\frac{1}{\gamma}\dot{u}(\nu)\right), \text{ for almost all } \nu \in \mathbb{R}^+. 
\end{split} \end{equation}
\noindent We define the function $z_\gamma$ by the relation $z_\gamma = x_\gamma \circ s_{1/\gamma}$ so that $$\dot{z}_\gamma (\nu) = \dot{x}_\gamma \big( s_{1/\gamma}(\nu) \big) \cdot \dot{s}_{1/\gamma}(\nu)=\dot{x}_\gamma (\gamma \nu) \cdot \gamma .$$
Substituting in (\ref{ecuacion1explicativa})  provides:
\begin{equation} \label{equation13} \begin{split}
&\dot{z}_\gamma\left(t\right)  =  -\sigma_{0}\frac{|\dot{u}(t)|}{g\left(\frac{\dot{u}(t)}{\gamma}\right)}z_\gamma(t)+\dot{u}(t), \text{ for almost all } t \in \mathbb{R}^+, \\
&z_\gamma(0)  =  x_{0}, \\
&y_\gamma(t)   =  \sigma_{0}z_\gamma(t)+\frac{\sigma_{1}}{\gamma}\dot{z}_\gamma(t)+f\left(\frac{\dot{u}(t)}{\gamma}\right), \text{ for almost all } t \in \mathbb{R}^+, 
\end{split}\end{equation}
where $y_\gamma=F_\gamma \circ s_{1/\gamma}$.

For a given $\gamma>0$, the corresponding  output-versus-input graph is the set $G_{u \circ s_\gamma} = \big \{ \big(u \circ s_\gamma(t),F_\gamma(t)\big),t \geq 0 \big \} = \big \{ \big(u(t),F_\gamma \circ s_{1/\gamma}(t)=y_\gamma(t)\big),t \geq 0 \big \}$. The hysteresis loop of system \eqref{equation13} is the output-versus-input graph obtained for very slow inputs (that is when $\gamma \rightarrow \infty$) in steady state. The next result, which is a straightforward combination of \cite[Proposition 5]{FI13} and \cite[Theorem 9]{Naser et al.(2015)}, describes the result of this convergence process. 

\begin{theorem}[\cite{FI13,Naser et al.(2015)}]\label{t:conv}
The following statements hold:
\begin{enumerate}
\item[\textnormal{(a)}] The sequence of functions $(y_\gamma )_{\gamma >0}$ converges in the space $L^\infty(\mathbb{R}^+,\mathbb{R})$ as $\gamma \rightarrow \infty$.  Denote $y^{\star}_{u}=\lim_{\gamma \rightarrow \infty} y_\gamma $, then for all $t \in \mathbb{R}^+$ we have
\begin{eqnarray}
y^{\star}_{u}(t) &=& \sigma_0 \mathrm{e}^{-\frac{\sigma_0}{g(0)}\rho_u(t)}\left(x_0 + \int_0^t \mathrm{e}^{\frac{\sigma_0}{g(0)}\rho_u(\tau)} \dot{u}(\tau)d\tau\right), \label{eqFstar1}\\
\rho_u(t) &=& \int_0^t |\dot{u}(\tau)|d\tau. \label{eqFstar2}
\end{eqnarray}
\item[\textnormal{(b)}] For any $k \in \mathbb{N}$ define the function $y^{\star}_{u,k}\in L^\infty\big([0,T],\mathbb{R} \big)$ by $y^{\star}_{u,k}(t)=y^{\star}_{u}(t+kT)$, for all $t \in [0,T]$. The sequence of functions $(y^{\star}_{u,k})_{k \in \mathbb{N}}$ converges in the space $L^\infty\big([0,T],\mathbb{R}\big)$ as $k \rightarrow \infty$.  Denote $y^{\circ}_u=\lim_{k \rightarrow \infty} y^{\star}_{u,k}$, then 
\begin{equation} \label{eqFcirc1}
y^{\circ}_u(t) =\sigma_0 \mathrm{e}^{-\frac{\sigma_0}{g(0)}\rho_u(t)}\left( \frac{y^{\circ}_u(0)}{\sigma_0}+ \int_0^t \mathrm{e}^{\frac{\sigma_0}{g(0)}\rho_u(\tau)} \dot{u}(\tau)d\tau\right) \mbox{ for all } t \in [0,T].
\end{equation}
Moreover, $y^{\circ}_u(T)=y^{\circ}_u(0)$.
\end{enumerate}
\end{theorem}

 Statement (a) implies  that  the graphs $G_{u \circ s_\gamma}$ converge in a sense  precised in \cite[Lemma 9]{I09}  to the graph $G^\star_u = \big \{ \big(u(t),y^{\star}_{u}(t)\big),t \geq 0 \big \}$ as $\gamma \rightarrow \infty$.  The  hysteresis loop is given by  the ``steady state'' of the parametrized curve $G^\star_u$ which by statement (b) is  the set
\begin{equation} \label{eqhystloopLuGre}
G_u^\circ = \big \{ \big(u(t),y^{\circ}_u(t)\big),t \in [0,T] \big \}.
\end{equation}

Moreover, Theorem \ref{t:conv} (b)  gives
\begin{eqnarray*}
y^{\circ}_u(T)&=&\sigma_0 \mathrm{e}^{-\frac{\sigma_0}{g(0)}\rho_u(T)}\left( \frac{y^{\circ}_u(0)}{\sigma_0}+ \int_0^T \mathrm{e}^{\frac{\sigma_0}{g(0)}\rho_u(\tau)} \dot{u}(\tau)d\tau\right),\\
y^{\circ}_u(T)&=&y^{\circ}_u(0),
\end{eqnarray*}
which leads to
\begin{equation} \label{eqy(0)circ}
y^{\circ}_u(0)=\frac{\sigma_0\mathrm{e}^{-\frac{\sigma_0}{g(0)}\rho_u(T)}\int_0^T \mathrm{e}^{\frac{\sigma_0}{g(0)}\rho_u(\tau)} \dot{u}(\tau)d\tau}{1-\mathrm{e}^{-\frac{\sigma_0}{g(0)}\rho_u(T)}}.
\end{equation}
Equations (\ref{eqFcirc1}) and (\ref{eqy(0)circ}) provide the analytical expression of the hysteresis loop (\ref{eqhystloopLuGre}).  This expression includes both the major loop and the minor loop. 

Observe that for the LuGre model neither $\sigma_1$ nor $f$ intervene in the equation of the hysteresis loop, and only the value $g(0)$ appears in this equation. Also note that Equations (\ref{eqFcirc1})--(\ref{eqy(0)circ}) are also valid for the Dahl model since the latter is a particular case of the LuGre model.

In Example 2  of Section \ref{ss:ex2} the reader can find a  detailed illustration of the concepts presented in this section.

\subsubsection{The normalized input}\label{ss:normalizedinput}

The hysteresis loop of the LuGre and the Dahl models is given in \eqref{eqhystloopLuGre}, and it is characterized by the function $y^{\circ}_u$ of Theorem \ref{t:conv} (b). Note that we are considering general input functions $u \in \mathbb{M}_{u_{\min,1},u_{\min,2},u_{\max,1}, u_{\max,2},t_1,t_2,t_3,t_4}$ that may not allow an explicit calculation of the integral present in Equation (\ref{eqFcirc1}).  To get an explicit calculation of that integral we follow the approach of \cite{I09} and \cite{Naser et al.(2015)}  that leads  to the explicit expression of the hysteresis loop by using 
the  so-called \emph{normalized input} $\psi_u$ associated to  $u$. 
The use of the normalized input will give another parametrization of the curve in \eqref{eqhystloopLuGre}, an explicit one.

\begin{figure}[H]
                \centering
                \includegraphics[scale=0.2]{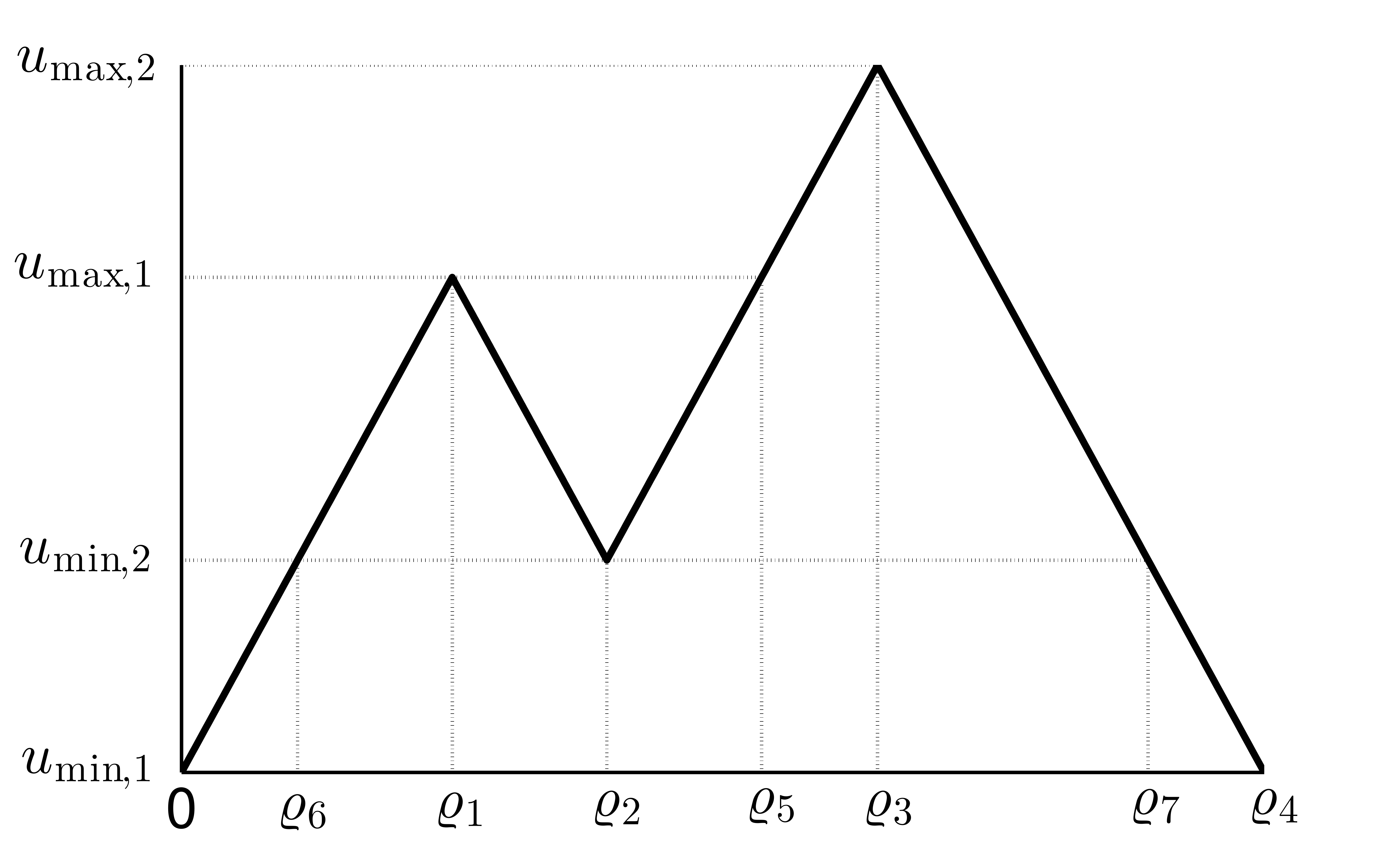}
                \caption{$\psi_u(\varrho)$ versus $\varrho$.}
                \label{MLX4}
\end{figure}

According to \cite{I09},  the normalized input associated to 
$u$ is a piecewise-linear function $\psi_u \in W^{1;\infty}(\mathbb{R}^+,\mathbb{R})$ that satisfies 
\begin{equation}\label{e:conjugacion}
\psi_u\big(\rho_u(t)\big)=u(t)\mbox{ for all } t\in\mathbb{R}^+, 
\end{equation} where
$
\rho_u(t) = \int_0^t |\dot{u}(\tau)|d\tau
$ is the \emph{variation function} of $u$. Note that $\rho_u$ is   strictly increasing so that it is invertible, and $\rho_u^{-1}$ is also strictly increasing. From equation (\ref{e:conjugacion}) it comes that $\psi_u = u \circ \rho_u^{-1}$ so that $\psi_u$ is strictly increasing on the interval $[0,\varrho_1]$, where $\varrho_1=\rho_u(t_1)$. Thus $\dot{\psi}_u(\varrho) \geq 0$ when $\varrho \in (0,\varrho_1)$  and $\dot{\psi}_u(\varrho)$ exists. On the other hand, by  \cite[Lemma 2]{I09}, the set on which $\dot{\psi}_u$ is not defined or is different from $\pm 1$ has measure zero. 
Thus $\dot{\psi}_u(\varrho) = 1$ for almost all $\varrho \in (0,\varrho_1)$. Using the fact that $\psi_u$ is absolutely continuous we obtain from the Fundamental Theorem of Calculus that $$ \psi_u(\varrho)-\psi_u(0)=\int_0^\varrho \dot{\psi}_u(\tau)\,\text{d}\tau=\int_0^\varrho \text{d}\tau=\varrho, \text{ for all } \varrho \in [0,\varrho_1].$$
Taking into account that $\psi_u\big(\rho_u(0)\big)=u(0)$ it comes that $\psi_u(0)=u_{\min,1}$ so that $$\psi_u(\varrho)=\varrho+u_{\min,1}, \text{ for all } \varrho \in [0,\varrho_1].$$
\noindent The details for the intervals $[\varrho_1,\varrho_2]$, $[\varrho_2,\varrho_3]$, and $[\varrho_3,\varrho_4]$ are given hereafter. 
\begin{itemize}
\item By definition of $u$ we have that $u$ is strictly increasing on the interval $[0,t_1]$ so that 
$$\varrho_1 = \rho_u(t_1) = \int_0^{t_1} |\dot{u}(t)| \text{d}t = \int_0^{t_1} \dot{u}(t) \text{d}t = u(t_1)-u(0) = u_{\max,1}-u_{\min,1}.$$
Also, from  $\psi_u(\varrho)=\varrho +u_{\min,1}$ for $\varrho \in [0,\varrho_1]$ we get $$\psi_u(\varrho_1)=\varrho_1+u_{\min,1}=u_{\max,1}.$$
Note that $\rho_u$ is   strictly increasing so that it is invertible, and $\rho_u^{-1}$ is also strictly increasing. From equation (\ref{e:conjugacion}) it comes that $\psi_u = u \circ \rho_u^{-1}$ so that $\psi_u$ is strictly decreasing on the interval $[\varrho_1,\varrho_2]$, where $\varrho_2=\rho_u(t_2)$. Thus $\dot{\psi}_u(\varrho) \leq 0$ when $\varrho \in (\varrho_1,\varrho_2)$  and $\dot{\psi}_u(\varrho)$ exists. On the other hand, by  \cite[Lemma 2]{I09}, the set on which $\dot{\psi}_u$ is not defined or is different from $\pm 1$ has measure zero. 
Thus $\dot{\psi}_u(\varrho) = -1$ for almost all $\varrho \in (\varrho_1,\varrho_2)$. Using the fact that $\psi_u$ is absolutely continuous we obtain from the Fundamental Theorem of Calculus that $$ \psi_u(\varrho)-\psi_u(\varrho_1)=\int_{\varrho_1}^\varrho \dot{\psi}_u(\tau)\,\text{d}\tau=\int_{\varrho_1}^\varrho -1 \;\text{d}\tau=\varrho_1-\varrho, \text{ for all } \varrho \in [\varrho_1,\varrho_2],$$
which leads to $$\psi_u(\varrho)=\psi_u(\varrho_1)+\varrho_1-\varrho=u_{\max,1}+\varrho_1-\varrho.$$
\item By definition of $u$ we have that $u$ is strictly decreasing on the interval $[t_1,t_2]$ so that 
\begin{eqnarray*}
\varrho_2 &=& \rho_u(t_2) = \int_0^{t_2} |\dot{u}(t)| \text{d}t = \underbrace{\int_0^{t_1} |\dot{u}(t)| \text{d}t}_{\varrho_1}+\int_{t_1}^{t_2} -\dot{u}(t) \text{d}t \\ &=&\varrho_1+ u(t_1)-u(t_2) = \varrho_1+u_{\max,1}-u_{\min,2}.
\end{eqnarray*}
Also, from  $\psi_u(\varrho)=u_{\max,1}+\varrho_1-\varrho$ for $\varrho \in [\varrho_1,\varrho_2]$ we get $$\psi_u(\varrho_2)=u_{\max,1}+\varrho_1-\varrho_2=u_{\min,2}.$$
 From equation (\ref{e:conjugacion}) it comes that $\psi_u = u \circ \rho_u^{-1}$ so that $\psi_u$ is strictly increasing on the interval $[\varrho_2,\varrho_3]$, where $\varrho_3=\rho_u(t_3)$. Thus $\dot{\psi}_u(\varrho) \geq 0$ when $\varrho \in (\varrho_2,\varrho_3)$  and $\dot{\psi}_u(\varrho)$ exists. On the other hand, by  \cite[Lemma 2]{I09}, the set on which $\dot{\psi}_u$ is not defined or is different from $\pm 1$ has measure zero. 
Thus $\dot{\psi}_u(\varrho) = 1$ for almost all $\varrho \in (\varrho_2,\varrho_3)$. Using the fact that $\psi_u$ is absolutely continuous we obtain from the Fundamental Theorem of Calculus that $$ \psi_u(\varrho)-\psi_u(\varrho_2)=\int_{\varrho_2}^\varrho \dot{\psi}_u(\tau)\,\text{d}\tau=\int_{\varrho_2}^\varrho \text{d}\tau=\varrho-\varrho_2, \text{ for all } \varrho \in [\varrho_2,\varrho_3],$$
which leads to $$\psi_u(\varrho)=\psi_u(\varrho_2)+\varrho-\varrho_2=u_{\min,2}+\varrho-\varrho_2.$$
\item By definition of $u$ we have that $u$ is strictly increasing on the interval $[t_2,t_3]$ so that 
\begin{eqnarray*}
\varrho_3 &=& \rho_u(t_3) = \int_0^{t_3} |\dot{u}(t)| \text{d}t = \underbrace{\int_0^{t_2} |\dot{u}(t)| \text{d}t}_{\varrho_2}+\int_{t_2}^{t_3} \dot{u}(t) \text{d}t \\ &=&\varrho_2+ u(t_3)-u(t_2) = \varrho_2+u_{\max,2}-u_{\min,2}.
\end{eqnarray*}
Also, from  $\psi_u(\varrho)=u_{\min,2}+\varrho-\varrho_2$ for $\varrho \in [\varrho_2,\varrho_3]$ we get $$\psi_u(\varrho_3)=u_{\min,2}+\varrho_3-\varrho_2=u_{\max,2}.$$
 From equation (\ref{e:conjugacion}) it comes that $\psi_u = u \circ \rho_u^{-1}$ so that $\psi_u$ is strictly decreasing on the interval $[\varrho_3,\varrho_4]$, where $\varrho_4=\rho_u(t_4)$. Thus $\dot{\psi}_u(\varrho) \leq 0$ when $\varrho \in (\varrho_3,\varrho_4)$  and $\dot{\psi}_u(\varrho)$ exists. On the other hand, by  \cite[Lemma 2]{I09}, the set on which $\dot{\psi}_u$ is not defined or is different from $\pm 1$ has measure zero. 
Thus $\dot{\psi}_u(\varrho) = -1$ for almost all $\varrho \in (\varrho_3,\varrho_4)$. Using the fact that $\psi_u$ is absolutely continuous we obtain from the Fundamental Theorem of Calculus that $$ \psi_u(\varrho)-\psi_u(\varrho_3)=\int_{\varrho_3}^\varrho \dot{\psi}_u(\tau)\,\text{d}\tau=\int_{\varrho_3}^\varrho -1 \;\text{d}\tau=\varrho_3-\varrho, \text{ for all } \varrho \in [\varrho_3,\varrho_4],$$
which leads to $$\psi_u(\varrho)=\psi_u(\varrho_3)+\varrho_3-\varrho=u_{\max,2}+\varrho_3-\varrho.$$
\end{itemize}
As a summary, we have
\begin{equation}\label{e:uuniversal}
\psi_u(\varrho)=\begin{cases}\begin{array}{ll}
\varrho+u_{\min,1}& \mbox{for } \varrho\in[0,\varrho_1],\\
-\varrho+\varrho_1+u_{\max,1}& \mbox{for } \rho\in[\varrho_1,\varrho_2],\\
\varrho-\varrho_2+u_{\min,2}& \mbox{for } \varrho \in[\varrho_2,\varrho_3],\\
-\varrho+\varrho_3+u_{\max,2}& \mbox{for } \varrho \in[\varrho_3,\varrho_4],
\end{array}
\end{cases}
\end{equation}
 where
$\varrho_1=\;u_{\max,1}-u_{\min,1}>0,$
$\varrho_2=\;\varrho_1+u_{\max,1}-u_{\min,2}>\varrho_1,$
$\varrho_3=\;\varrho_2+u_{\max,2}-u_{\min,2}>\varrho_2,$ and
$\varrho_4=\;\varrho_3+u_{\max,2}-u_{\min,1}>\varrho_3$. The function $\psi_u$ is continuous and $\varrho_4$--periodic. Its graph in the interval $[0,\varrho_4] $ is displayed in Figure \ref{MLX4}.

\subsubsection{Analytic expression of the hysteresis and minor loops} \label{Analytical expression of the hysteresis and minor loops}

Applying Theorem \ref{t:conv} (b) (Equation \eqref{eqFcirc1}) to the particular input $\psi_u$, and denoting for simplicity $y^\circ:=y^{\circ}_{\psi_u}$ we obtain that
$$ y^{\circ}(\varrho)  =\sigma_0 \mathrm{e}^{-\frac{\sigma_0}{g(0)}\varrho}\left( \frac{y^\circ(0)}{\sigma_0}+ \int_0^\varrho \mathrm{e}^{\frac{\sigma_0}{g(0)}\tau} \dot{\psi}_u(\tau)d\tau\right) \mbox{ for  } \varrho \in [0,\varrho_4].
$$
Since this expression can be explicitly integrated, we obtain
\begin{eqnarray}
y^\circ(\varrho)&=& \mathrm{e}^{-\frac{\sigma_0}{g(0)}\varrho}\left( y^\circ(0)+ g(0)\left[ \mathrm{e}^{\frac{\sigma_0}{g(0)}\varrho} -1\right]\right) \mbox{ for } \varrho \in [0,\varrho_1], \mbox{ with } \label{ycircdezeroML}\\
y^\circ(0) &=&g(0)\,\frac{\mathrm{e}^{-\frac{\sigma_0}{g(0)}\varrho_4}}{1-\mathrm{e}^{-\frac{\sigma_0}{g(0)}\varrho_4}}\left(  2\mathrm{e}^{\frac{\sigma_0}{g(0)}\varrho_1} -2\mathrm{e}^{\frac{\sigma_0}{g(0)}\varrho_2}+2\mathrm{e}^{\frac{\sigma_0}{g(0)}\varrho_3}-\mathrm{e}^{\frac{\sigma_0}{g(0)}\varrho_4}-1\right); \notag
\end{eqnarray}
\begin{eqnarray}
y^\circ(\varrho)&=&\mathrm{e}^{-\frac{\sigma_0}{g(0)}(\varrho-\varrho_1)}\left( y^\circ(\varrho_1)- g(0)\left[ \mathrm{e}^{\frac{\sigma_0}{g(0)}(\varrho-\varrho_1)} -1\right]\right) \mbox{ for } \varrho \in [\varrho_1,\varrho_2],  \mbox{ with }\label{[rho1rho2ycirc]}\\
y^\circ(\varrho_1)&=&\mathrm{e}^{-\frac{\sigma_0}{g(0)}\varrho_1}\left( y^\circ(0)+ g(0)\left[ \mathrm{e}^{\frac{\sigma_0}{g(0)}\varrho_1} -1\right]\right); \notag 
\end{eqnarray}
\begin{eqnarray}
y^\circ(\varrho)&=&\mathrm{e}^{-\frac{\sigma_0}{g(0)}(\varrho-\varrho_2)}\left( y^\circ(\varrho_2)+ g(0)\left[ \mathrm{e}^{\frac{\sigma_0}{g(0)}(\varrho-\varrho_2)} -1\right]\right) \mbox{ for } \varrho \in [\varrho_2,\varrho_3], \mbox{ with } \label{[rho2rho3ycirc]}\\
y^\circ(\varrho_2)&=&\mathrm{e}^{-\frac{\sigma_0}{g(0)}(\varrho_2-\varrho_1)}\left( y^\circ(\varrho_1)- g(0)\left[ \mathrm{e}^{\frac{\sigma_0}{g(0)}(\varrho_2-\varrho_1)} -1\right]\right);\notag 
\end{eqnarray}
and
\begin{eqnarray}
y^\circ(\varrho)&=&\mathrm{e}^{-\frac{\sigma_0}{g(0)}(\varrho-\varrho_3)}\left( y^\circ(\varrho_3)- g(0)\left[ \mathrm{e}^{\frac{\sigma_0}{g(0)}(\varrho-\varrho_3)} -1\right]\right) \mbox{ for } \varrho \in [\varrho_3,\varrho_4], \mbox{ with } \label{[rho3rho4ycirc]}\\
y^\circ(\varrho_3)&=&\mathrm{e}^{-\frac{\sigma_0}{g(0)}(\varrho_3-\varrho_2)}\left( y^\circ(\varrho_2)+ g(0)\left[ \mathrm{e}^{\frac{\sigma_0}{g(0)}(\varrho_3-\varrho_2)} -1\right]\right). \notag 
\end{eqnarray}

Finally, observe that
the hysteresis loop of the LuGre model that corresponds to the input $\psi_u$ is the set 
\begin{equation} \label{analytMLeq}
G_{\psi_u}^\circ = \big \{ \big(\psi_u(\varrho),y^\circ(\varrho)\big),\,\varrho \in [0,\varrho_4] \big \}.
\end{equation}
Taking into account the fact that $y^\circ_u=y^\circ \circ \rho_u$ and $u=\psi_u \circ \rho_u$ it comes from \cite[Lemma 8]{I09}  that $G_{\psi_u}^\circ=G_u^\circ = \big \{ \big(u(t),y_u^\circ(t)\big),\,t \in [0,T] \big \}$, thus proving statement (a) of Theorem \ref{t:main}.

To prove statement (b) of Theorem \ref{t:main} observe that  the minor loop corresponding to the input $\psi_u$ is the part of the hysteresis loop \eqref{analytMLeq} that corresponds to $\varrho \in [\varrho_1,\varrho_5]$, where 
$$
\varrho_5=\rho_u(t_5)=\int_{0}^{t_5} |\dot{u}(t)|dt=\varrho_2+u_{\max,1}-u_{\min,2}\in(\varrho_2,\varrho_3),
$$
and where $t_2<t_5<t_3$ is the time such that $u(t_5)=u(t_1)$, see Figure \ref{f:figugen}. 
This set is the union of the two arcs $\big\{\big(\psi_u(\varrho),y^\circ(\varrho)\big),\, \varrho \in [\varrho_1,\varrho_2]\big\}$ and $\big\{\big(\psi_u(\varrho),y^\circ(\varrho)\big),\,\varrho \in [\varrho_2,\varrho_5]\big\}$.

With this argument, the proof of Theorem \ref{t:main} is complete. \\

We remark that the explicit construction of the hysteretic loop, and therefore the identification of the arcs corresponding to the minor loops given in the proof of Theorem \ref{t:main} can be generalized to  
multimodal input functions giving rise to hysteresis loops with many minor loops. This can be done using the normalized input and Equation \eqref{eqFcirc1}.

\section{Examples}\label{s:numsim}

\subsection{Example 1: an approach to the concept of hysteresis loop}\label{ss:ex1}

A hysteresis system is  one for which the  output-versus-input graph presents a loop in the steady state for slow inputs \cite{Ikhouane2017}. The way to obtain the hysteresis loop that corresponds to a given input is as follows. Consider a periodic input $t \rightarrow u(t)$. Composing this input with the time-scale change $t \rightarrow t/\gamma$ provides a new input $u_\gamma:t \rightarrow u( t/\gamma)$. This new input gives rise to an output $y_\gamma(t)$  such that the corresponding output-versus-input graph  $\big\{\big(u_\gamma(t),y_\gamma(t)\big),t \geq 0\big\}$ converges to a fixed curve -the hysteresis loop- in steady state when $\gamma \rightarrow \infty$. 
Our aim in this section is to illustrate this process using an example. \\

\begin{figure}[H]
                \centering
                \includegraphics[scale=0.15]{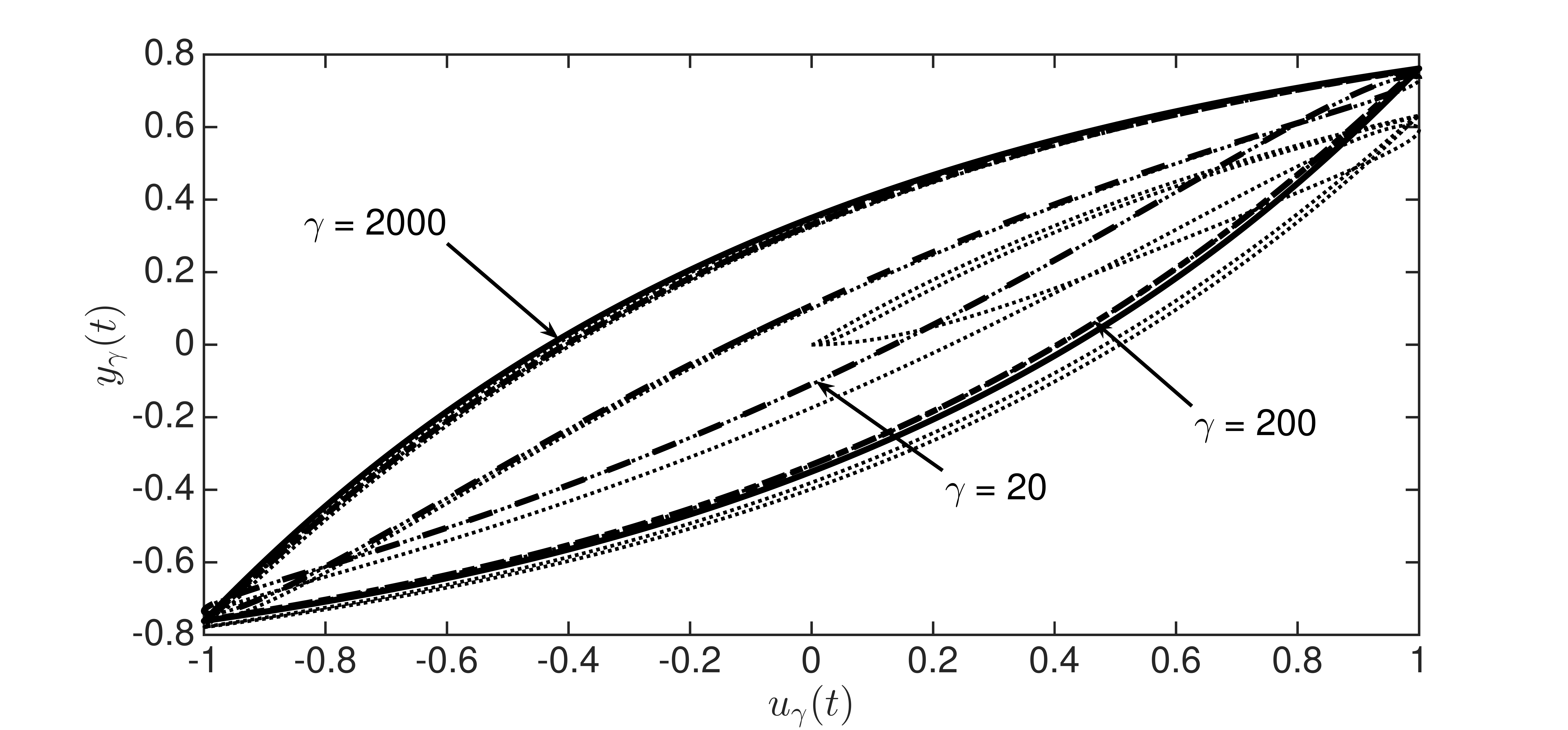}
                \caption{Output of system (\ref{eqPCCp}), $y_\gamma$ versus $u_\gamma$. Dotted: transient; solid: steady state for $\gamma = 2000$; dashed: steady state for $\gamma = 200$; dash-dotted: steady state for $\gamma = 20$.}
                \label{Figguraseis}
\end{figure}

Consider for instance the following system constructed using the Dahl model:
\begin{equation}\label{eqPCCp}\begin{split}
\dot{x}_\gamma(t)=&\;\dot{u}_\gamma(t)-|\dot{u}_\gamma(t)|x_\gamma(t),\\
\dot{y}_\gamma(t)=&-y_\gamma(t)+x_\gamma(t),\\
x_\gamma(0)=& \;0,\quad y_\gamma(0)=0, 
\end{split}\end{equation}
with input $u_\gamma(t)=\sin(2\pi t/\gamma)$  and output $y_\gamma(t)$. Figure \ref{Figguraseis} provides the output-versus-input graph  $\big\{\big(u_\gamma(t),y_\gamma(t)\big),\, t \geq 0\big\}$ of system (\ref{eqPCCp}) for increasing values of $\gamma$.
 \noindent It can be seen that, as $\gamma \rightarrow \infty$, the steady-state part of the output-versus-input graph
converges to a fixed closed curve. This curve is  the hysteresis loop of system (\ref{eqPCCp}).

 \begin{figure}[H]
                \centering
                \includegraphics[scale=0.15]{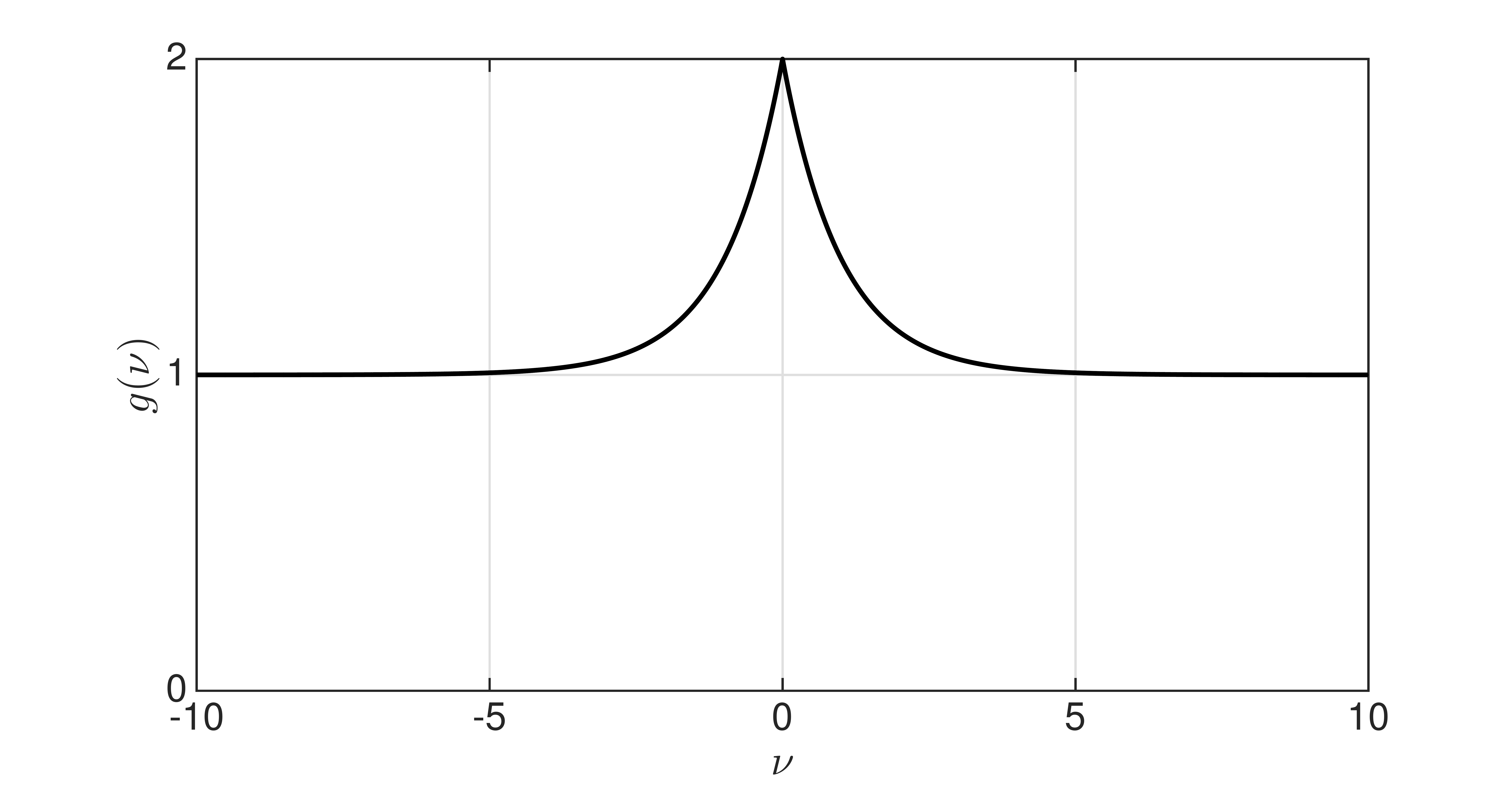}
                \caption{ The macrodamping friction function $g(\nu)$ versus $\nu$.}
                \label{mlfig2}
\end{figure}

\subsection{Example 2: the hysteresis loop of the LuGre model}\label{ss:ex2}

 The aim of this section is to illustrate the concepts presented in Section \ref{ss:hystloop} by means of numerical simulations.

 Following \cite{JC08}, to  approximate the Stribeck effect  we set: 
$$
g(\nu)=F_{c}+\left(F_{s}-F_{c}\right)\mathrm{e}^{-\left|\nu/v_{s}\right|^{\beta}}\,\mbox{ for }\nu\in\mathbb{R},
$$
where $F_{c}>0$ is the Coulomb friction force, $F_{s}>0$ is the stiction force, $v_{s}>0$ is the Stribeck velocity, and $\beta$ is a strictly positive constant. The function $f$ is taken to be zero. The values of the different constants are taken to be  $\sigma_{0}=1$, $\sigma_1=1$,  $F_c=1$, $F_s=2$, $v_{s}=1$,  $\beta=1$, see Figure \ref{mlfig2}.

\begin{figure}[H]
                \centering
                \includegraphics[scale=0.25]{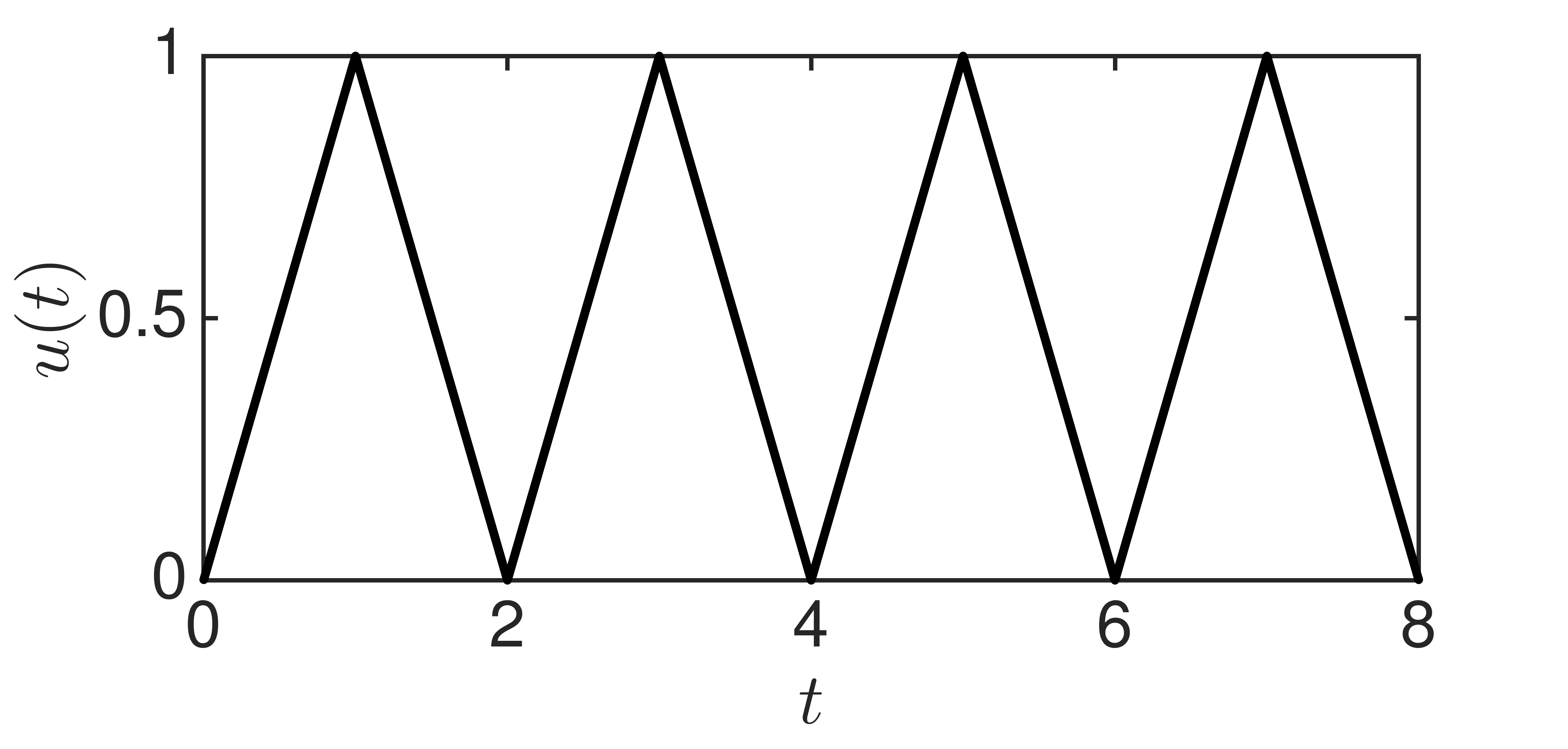}
                \caption{$u(t)$ versus $t$.}
                \label{mlfig1}
\end{figure}

The input is the continuous $2$-periodic piecewise-linear function defined by $u(t)=t$ for $t\in [0,1]$ and $u(t)=2-t$ for $t \in [1,2]$; see Figure \ref{mlfig1}. Observe that $\psi_u=u$.

\begin{figure}[H]
                \centering
                \includegraphics[scale=0.15]{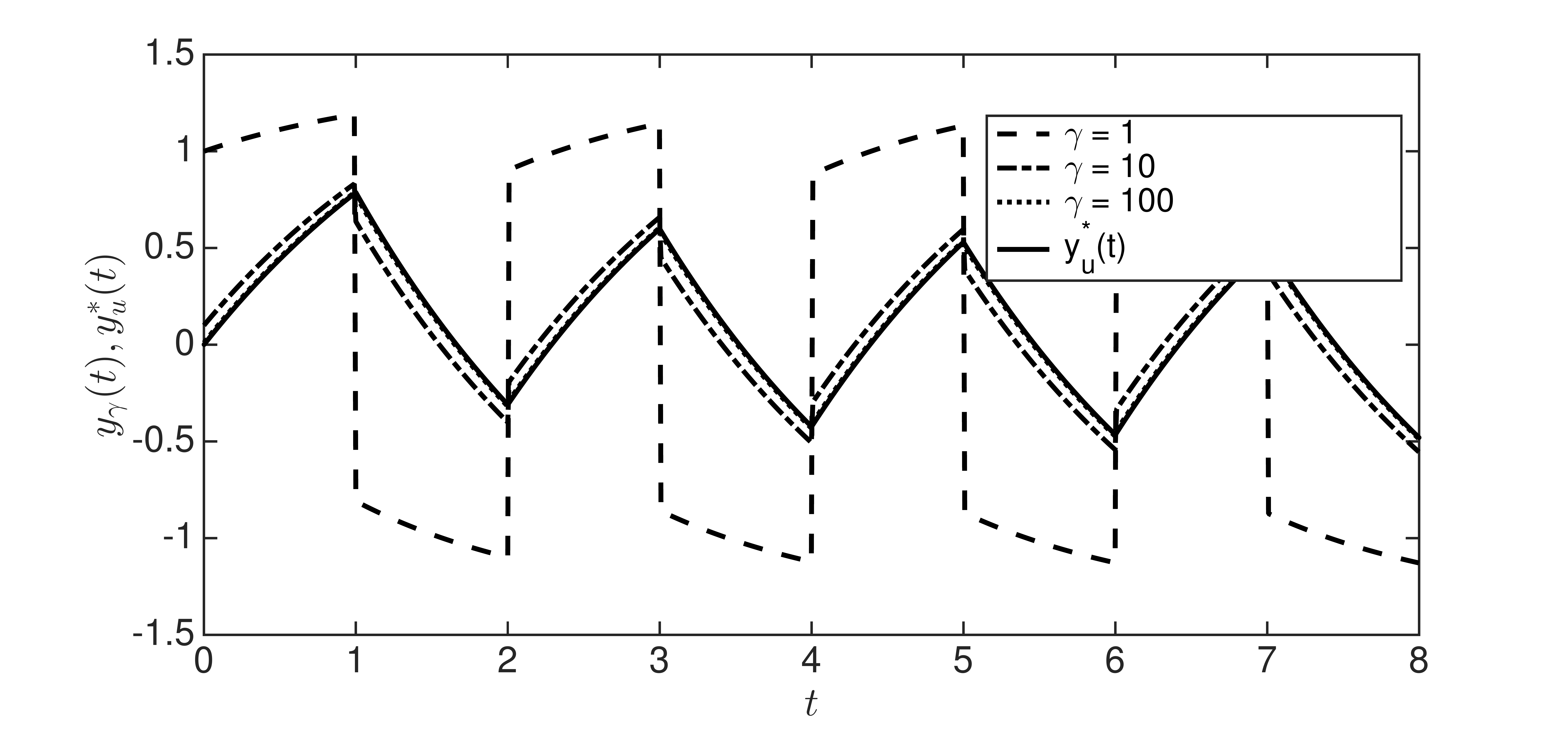}
                \caption{$y_\gamma(t)$ versus $t$. Dashed $\gamma=1$, dash-dotted $\gamma=10$, dotted $\gamma=100$; solid $y^{\star}_{u}(t)$ versus $t$.}
                \label{mlfig3}
\end{figure}

We take  $x\left(0\right)=x_0=0$. With these values we obtain $y_\gamma$ by a numerical integration of  Equations (\ref{equation13}) using Matlab solver \texttt{ode23s}.. Also, using Equations (\ref{eqFstar1})--(\ref{eqFstar2}) we obtain $y^{\star}_{u}$. Figure \ref{mlfig3} provides the plots $y_\gamma(t)$ versus $t$ for $\gamma=1,10,100$ along with the plot $y^{\star}_{u}(t)$ versus $t$. It can be seen that as $\gamma$ increases, $y_\gamma$ converges to $y^{\star}_{u}$.

\begin{figure}[H]
                \centering
                \includegraphics[scale=0.15]{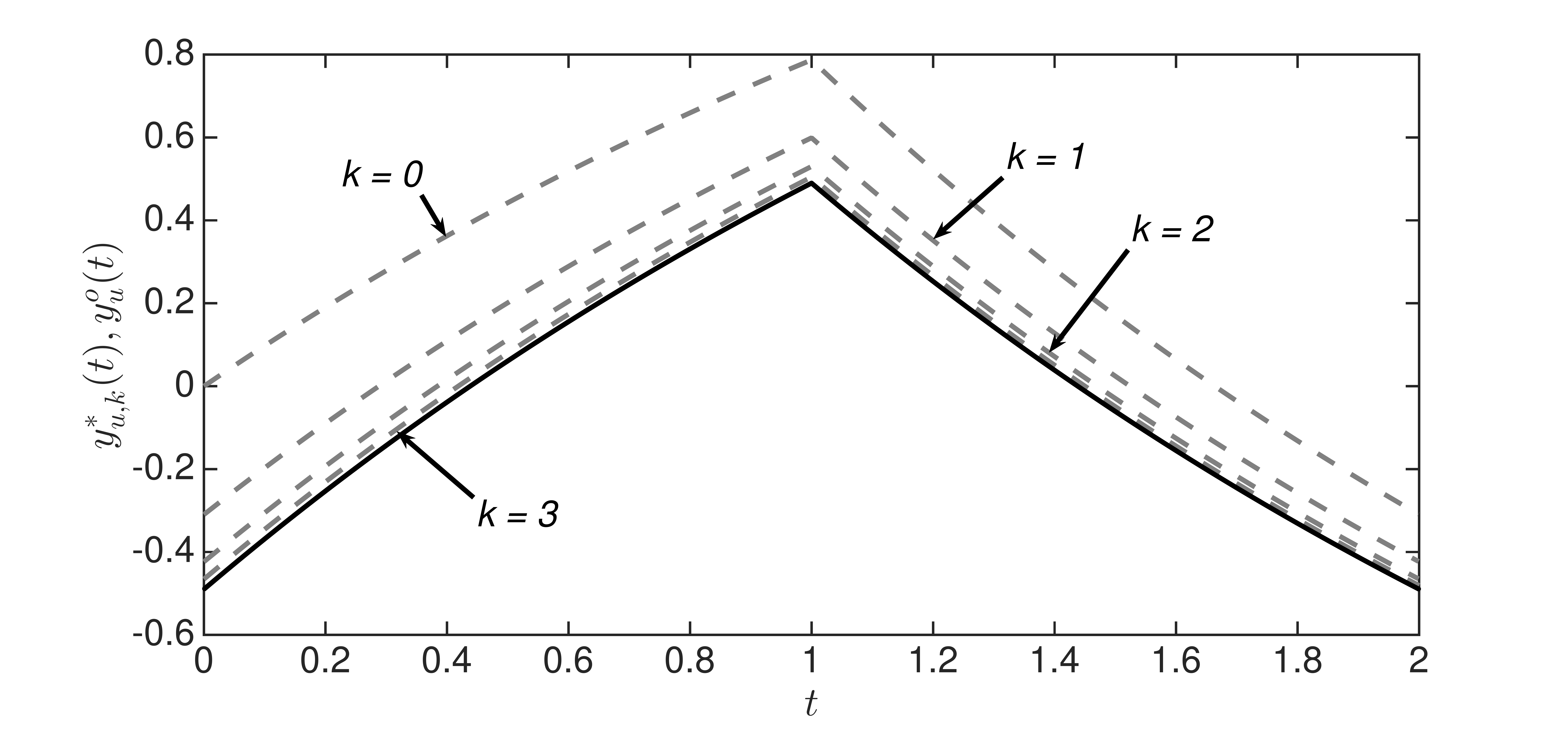}
                \caption{Convergence of the functions $y^{\star}_{u,k}$ to $y_u^\circ$. Dashed $y^{\star}_{u,k}(t)$ versus $t$ for $k=0,1,2,3$;  solid $y_u^\circ(t)$ versus $t$.}
                \label{mlfig4}
\end{figure}

The functions $y^{\star}_{u,k}$ are given by $y^{\star}_{u,k}(t)=y^{\star}_{u}(t+kT),t \in [0,T],k\in \mathbb{N}$ whilst $y_u^\circ$  is calculated from Equations (\ref{eqFcirc1}) and (\ref{eqy(0)circ}). Figure \ref{mlfig4} provides the plots $y^{\star}_{u,k}(t)$ versus $t$ for increasing values of $k$. It can be seen that $y^{\star}_{u,k}$ converges to $y_u^\circ$ as $k \rightarrow \infty$.

\begin{figure}[H]
                \centering
                \includegraphics[scale=0.15]{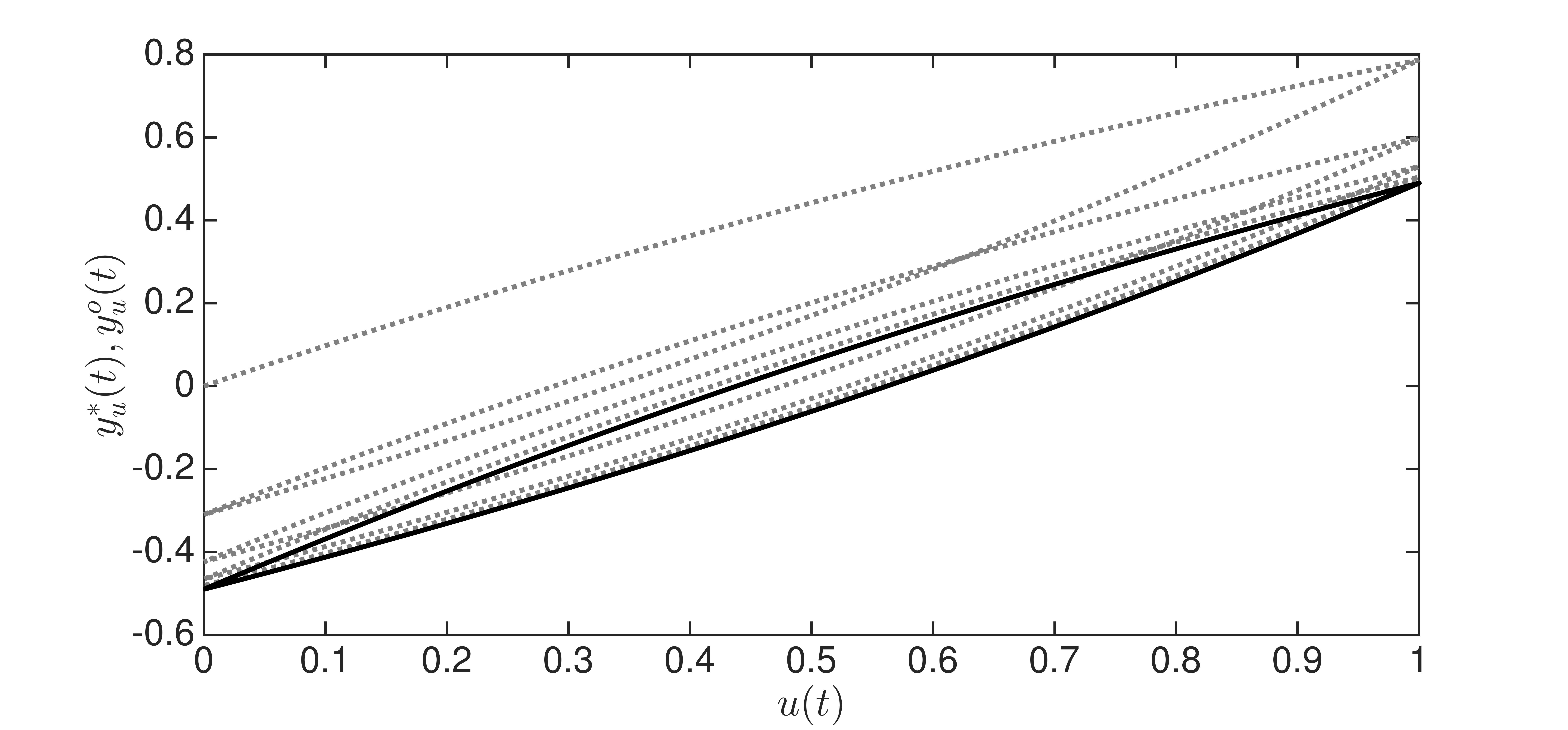}
                \caption{Convergence of the output-versus-input graph $G^\star_u$ to the hysteresis loop $G^\circ_u$.
                Dotted $G^\star_u$; solid $G^\circ_u$.}
                \label{mlfig5}
\end{figure} 

As in Example 1, the hysteresis loop is the output-versus-input graph  obtained for very slow inputs (that is when $\gamma \rightarrow \infty$) in steady state (that is when $k \rightarrow \infty$). For a given $\gamma$, the corresponding  output-versus-input graph is the set $G_{u \circ s_\gamma} = \big \{ \big(u \circ s_\gamma(t),F_\gamma(t)\big),t \geq 0 \big \} = \big \{ \big(u(t),F_\gamma \circ s_{1/\gamma}(t)=y_\gamma(t)\big),t \geq 0 \big \}$. Owing to Theorem \ref{t:conv} (a) and to \cite[Lemma 9]{I09} it comes that the graphs $G_{u \circ s_\gamma}$ converge in a sense detailed in \cite{I09} to the graph $G^\star_u = \big \{ \big(u(t),y^{\star}_{u}(t)\big),t \geq 0 \big \}$ as $\gamma \rightarrow \infty$.   Equations (\ref{eqFcirc1}) and (\ref{eqy(0)circ}) provide the analytic expression of the hysteresis loop (\ref{eqhystloopLuGre}).   Finally, Figure \ref{mlfig5} provides the graph $G^\star_u$ along with the hysteresis loop $G^\circ_u$.

\subsection{Example 3: process of convergence that leads to the hysteresis and minor loops} \label{ML0Numerical simulations}

The aim of Sections \ref{ML0Numerical simulations} and \ref{MLNumerical simulations} is to illustrate Theorem \ref{t:main} by means of numerical simulations. Section \ref{ML0Numerical simulations} focuses on the process of convergence that leads to the hysteresis loop. Section \ref{MLNumerical simulations} focuses on the variation of the  minor loop with the model's parameters.
 
 The function $g$ that characterizes the Stribeck effect is the same as in Section \ref{ss:ex2}. Also, the function $f$ is taken to be zero. 

\begin{figure}[H]
                \centering
                 \includegraphics[scale=0.15]{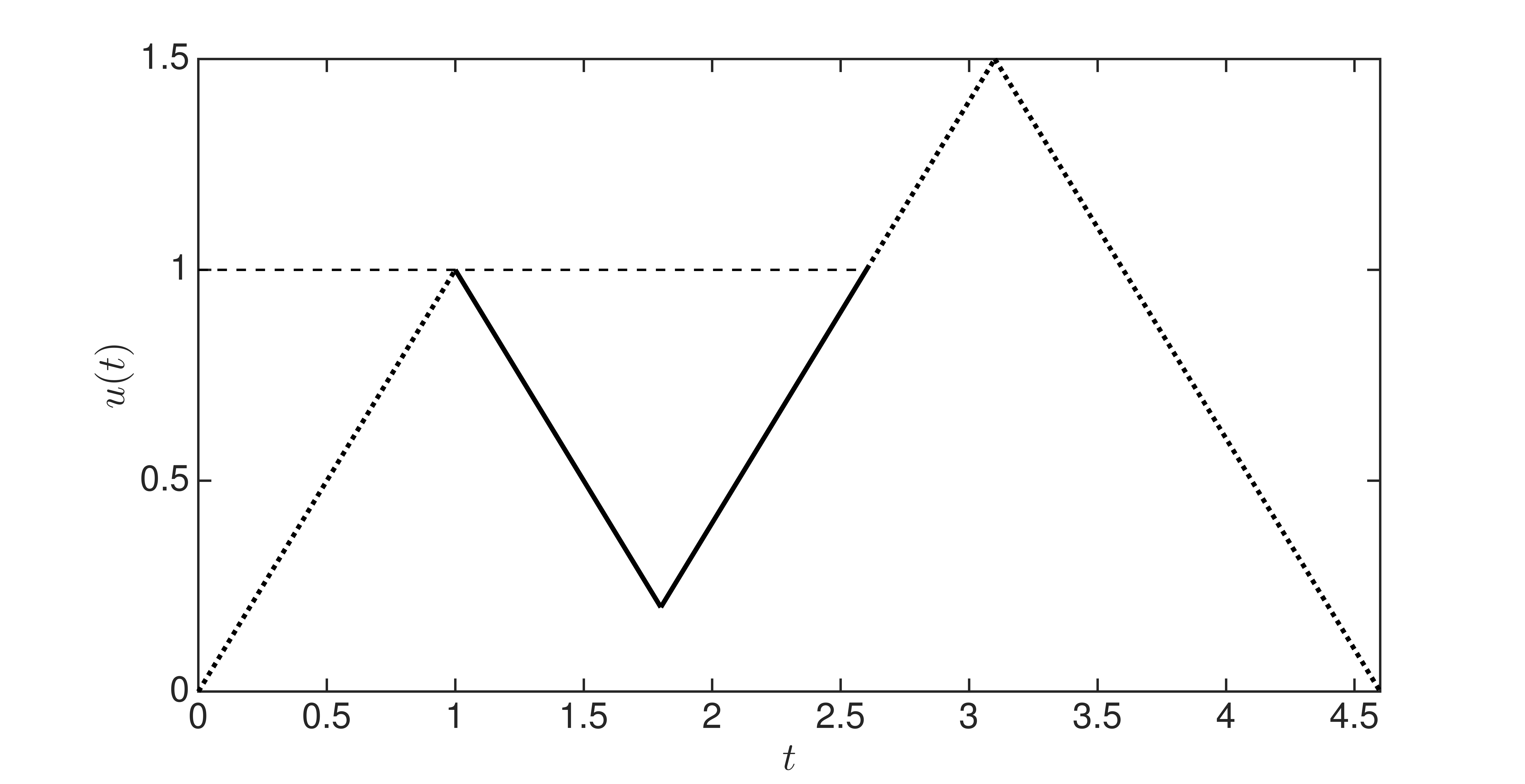}
                \caption{$u(t)$ versus $t$.}
                \label{figv2_ref4}
\end{figure}

The input is the continuous $\rho_4$-periodic piecewise-linear function defined by Equation (\ref{e:uuniversal}) where  $u_{\min,1}=0$, $u_{\min,2}=0.2$, $u_{\max,1}=1$,  $u_{\max,2}=1.5$, $\varrho_1=1$, $\varrho_2 =  1.8$, $\varrho_3 = 3.1$, $\varrho_4 = 4.6$, see Figure \ref{figv2_ref4}. Observe that $\psi_u=u$ and that the normalized variable $\varrho$ is equal to time $t$.

\begin{figure}[H]
                \centering
                 \includegraphics[scale=0.15]{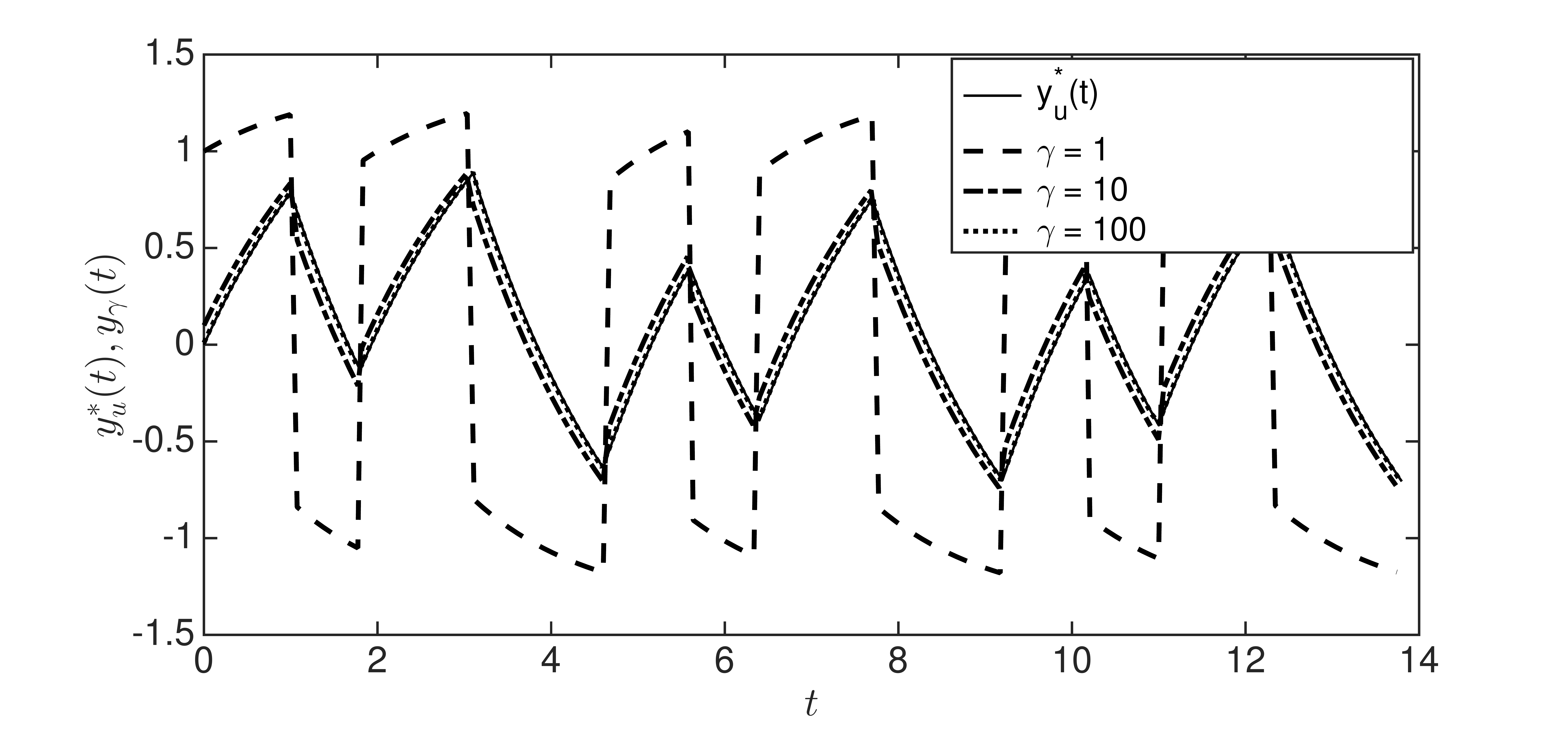}
                \caption{$y_\gamma(t)$ versus $t$. Dashed $\gamma=1$, dash-dotted $\gamma=10$, dotted $\gamma=100$; solid $y^{\star}_{u}(t)$ versus $t$.}
                \label{figv2_ref1}
\end{figure}

We take  $x\left(0\right)=x_0=0$. With these values we obtain $y_\gamma$ by a numerical integration of  Equations (\ref{equation13}) using Matlab solver \texttt{ode23s}. Also, using Equations (\ref{eqFstar1})--(\ref{eqFstar2}) we obtain $y^{\star}_{u}$. Figure \ref{figv2_ref1} provides the plots $y_\gamma(t)$ versus $t$ for $\gamma=1,10,100$ along with the plot $y^{\star}_{u}(t)$ versus $t$. It can be seen that as $\gamma$ increases, $y_\gamma$ converges to $y^{\star}_{u}$.

\begin{figure}[H]
                \centering
                 \includegraphics[scale=0.15]{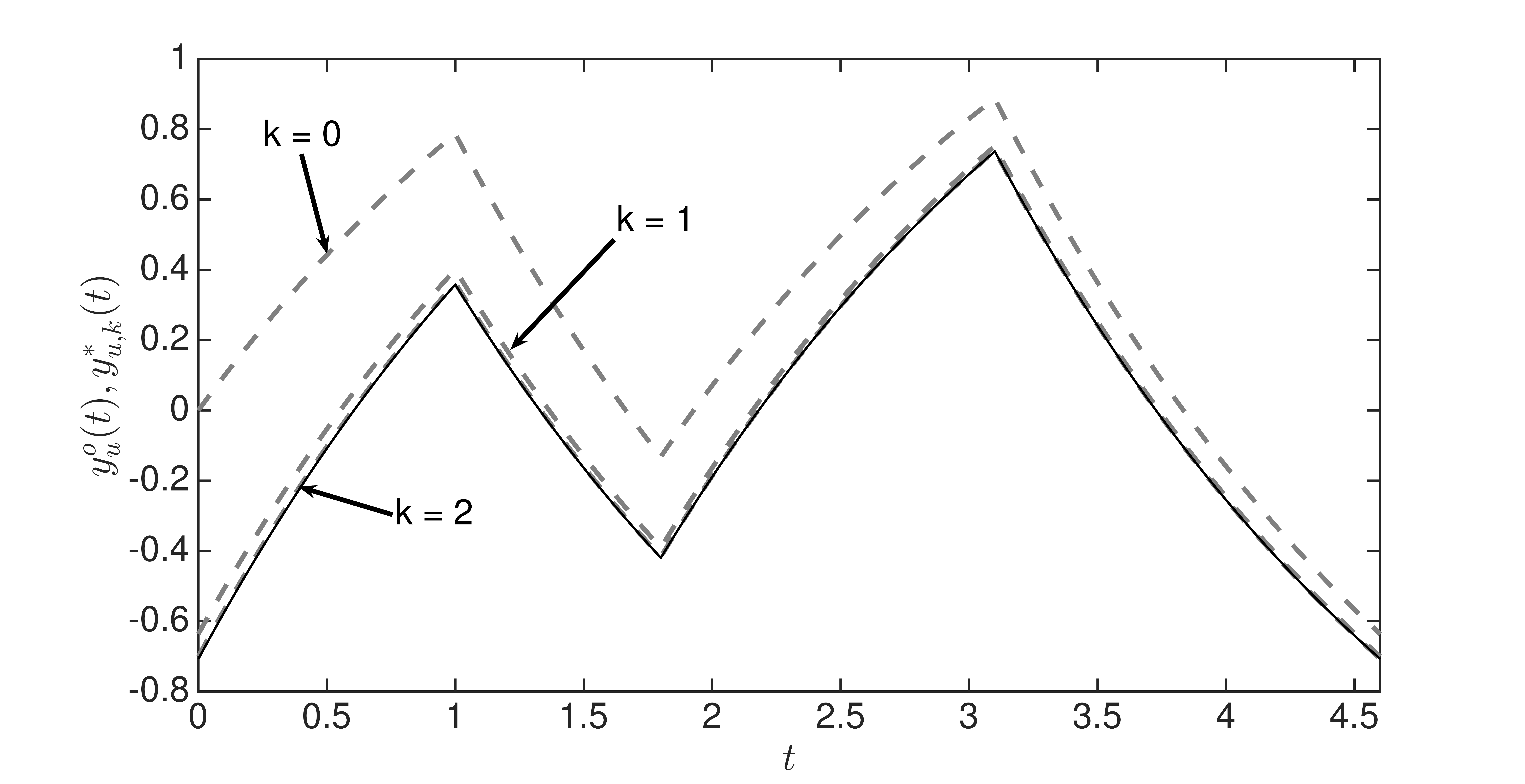}
                \caption{Convergence of the functions $y^{\star}_{u,k}$ to $y_u^\circ$. Dashed $y^{\star}_{u,k}(t)$ versus $t$ for $k=0,1,2$;  solid $y_u^\circ(t)$ versus $t$.}
                \label{figv2_ref2}
\end{figure}

The functions $y^{\star}_{u,k}$ are given by $y^{\star}_{u,k}(t)=y^{\star}_{u}(t+k\varrho_4),t \in [0,\varrho_4],k\in \mathbb{N}$ whilst $y_u^\circ$  is calculated from Equations (\ref{eqFcirc1}) and (\ref{eqy(0)circ}). Figure \ref{figv2_ref2} provides the plots $y^{\star}_{u,k}(t)$ versus $t$ for increasing values of $k$. It can be seen that $y^{\star}_{u,k}$ converges to $y_u^\circ$ as $k \rightarrow \infty$.

\begin{figure}[H]
                \centering
                 \includegraphics[scale=0.15]{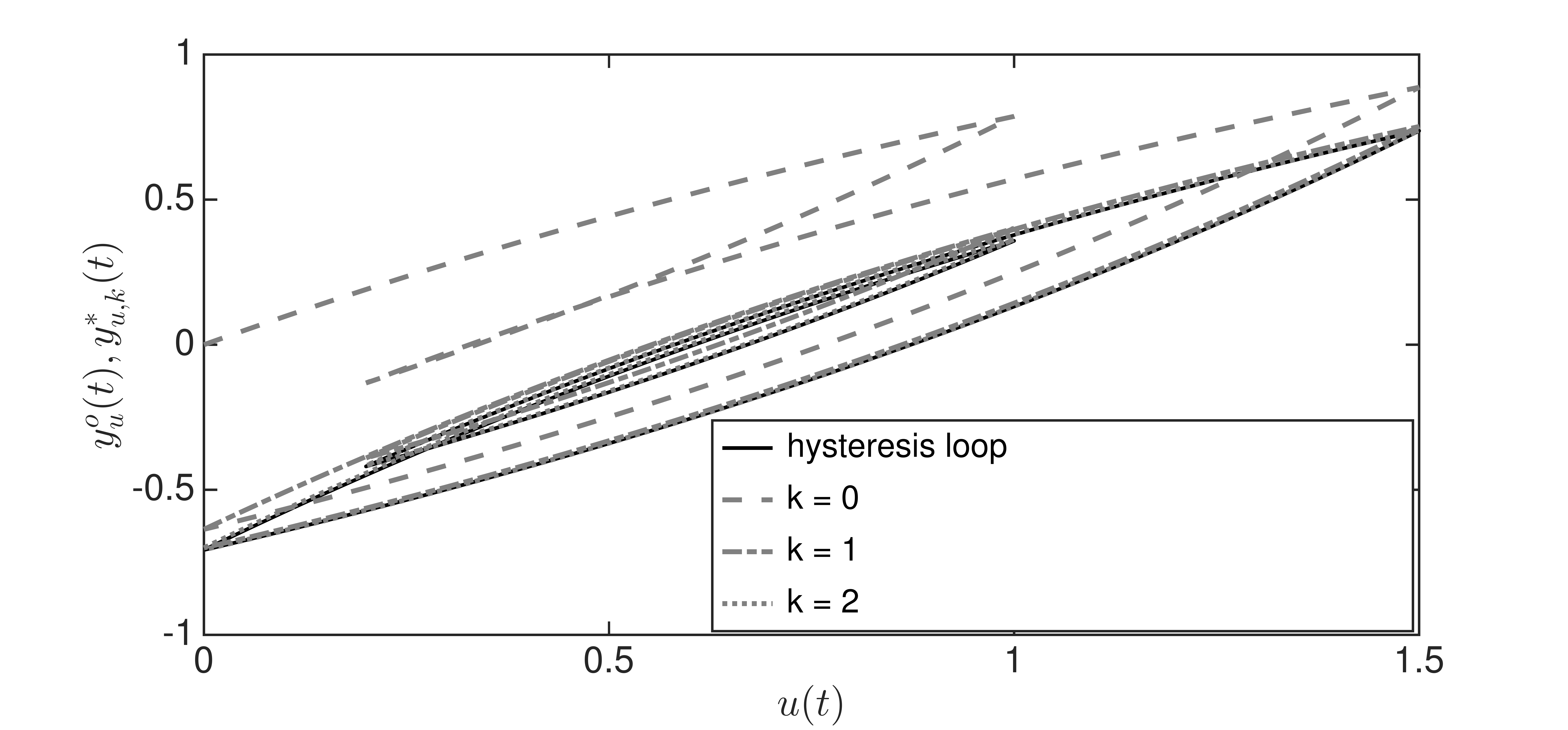}
                \caption{Convergence of the graphs $\big\{\big(u(t), y^{\star}_{u,k}(t)\big),t \in [0,\varrho_4]\big\}$ to the hysteresis loop $\big\{\big(u(t), y_u^{\circ}(t)\big),t \in [0,\varrho_4]\big\}$.}
                \label{figv2_ref3}
\end{figure}

Figure \ref{figv2_ref3} provides the graphs $\big\{\big(u(t), y^{\star}_{u,k}(t)\big),t \in [0,\varrho_4]\big\}$ for $k=0,1,2$. It can be seen that these graphs converge to the hysteresis loop $\big\{\big(u(t), y_u^{\circ}(t)\big),t \in [0,\varrho_4]\big\}$.

\subsection{Example 4: Variation of the minor loop with the model's parameters} \label{MLNumerical simulations}

\begin{figure}[H]
                \centering
                \includegraphics[scale=0.25]{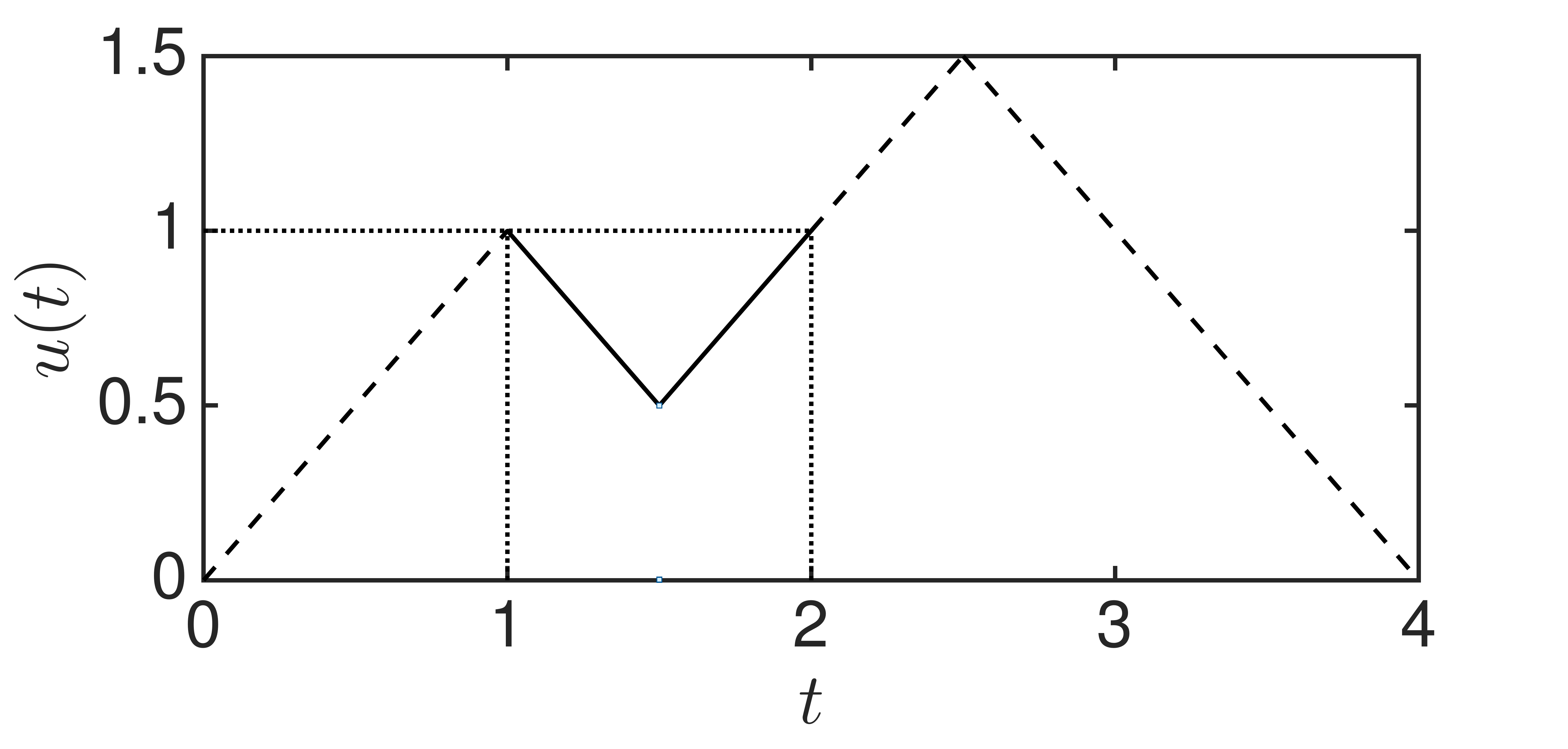}
                \caption{$u(t)$ versus $t$}
                \label{mlfig6}
\end{figure}

 We consider the LuGre model of  Section \ref{ss:ex2} with the value $\sigma_0=6$. The input $u$ is the one given in \eqref{e:uuniversal}
(thus a normalized one) with $u_{\min,1}=0$, $u_{\min,2}=0.5$, $u_{\max,1}=1$, $u_{\max,2}=1.5$,  with its corresponding values of $\varrho_i=t_i$   for $i=1,\ldots,5$,  see Figure \ref{mlfig6}.

\begin{figure}[H]
                \centering
                \includegraphics[scale=0.25]{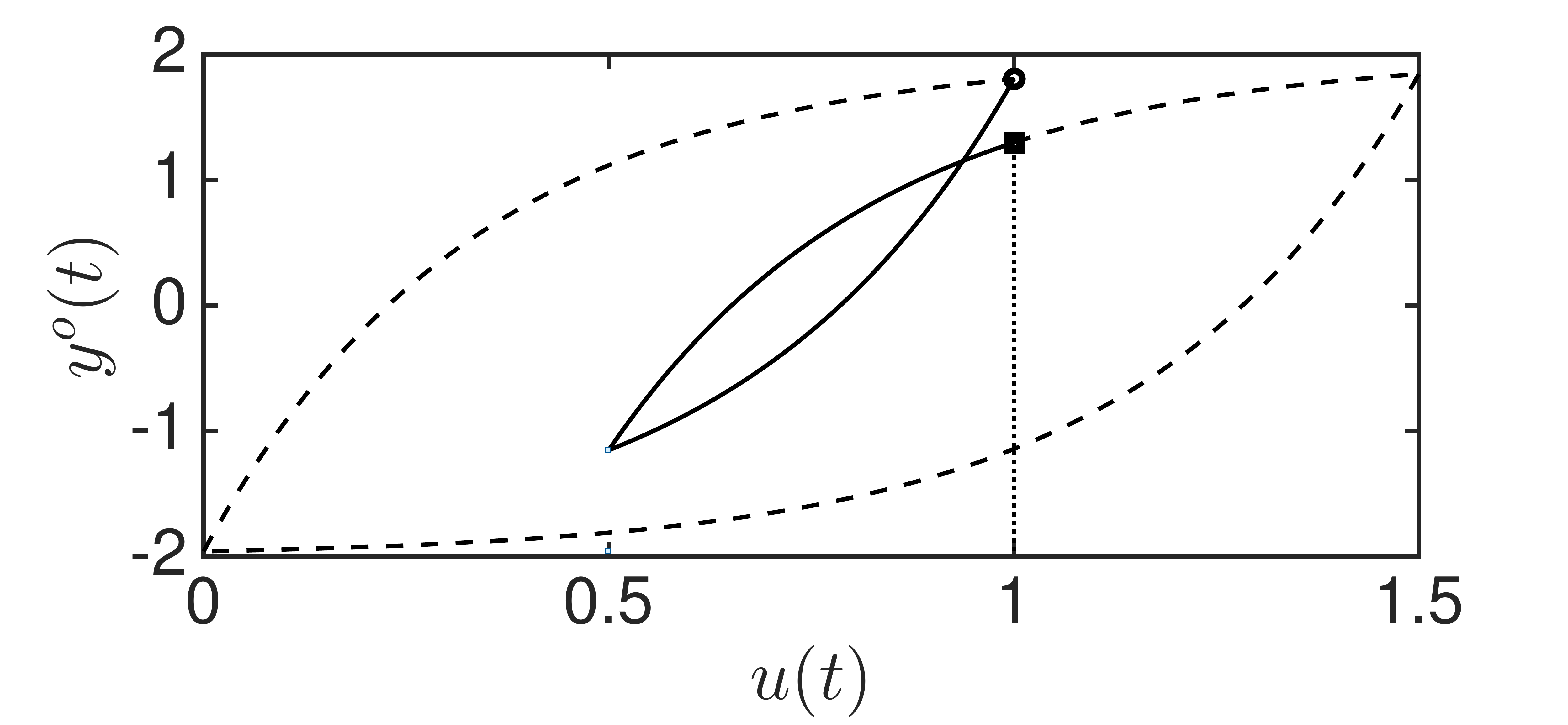}
                \caption{Hysteresis loop $G_u^\circ$. The marker \textit{open circle} corresponds to the point $\big(u(t_1),y^\circ(t_1)\big)$ whilst the marker \textit{rectangle} corresponds to the point $\big(u(t_5),y^\circ(t_5)\big)$. The minor loop is plotted in solid line.}
                \label{mlfig7}
\end{figure}

\begin{figure}[H]
                \centering
                \includegraphics[scale=0.25]{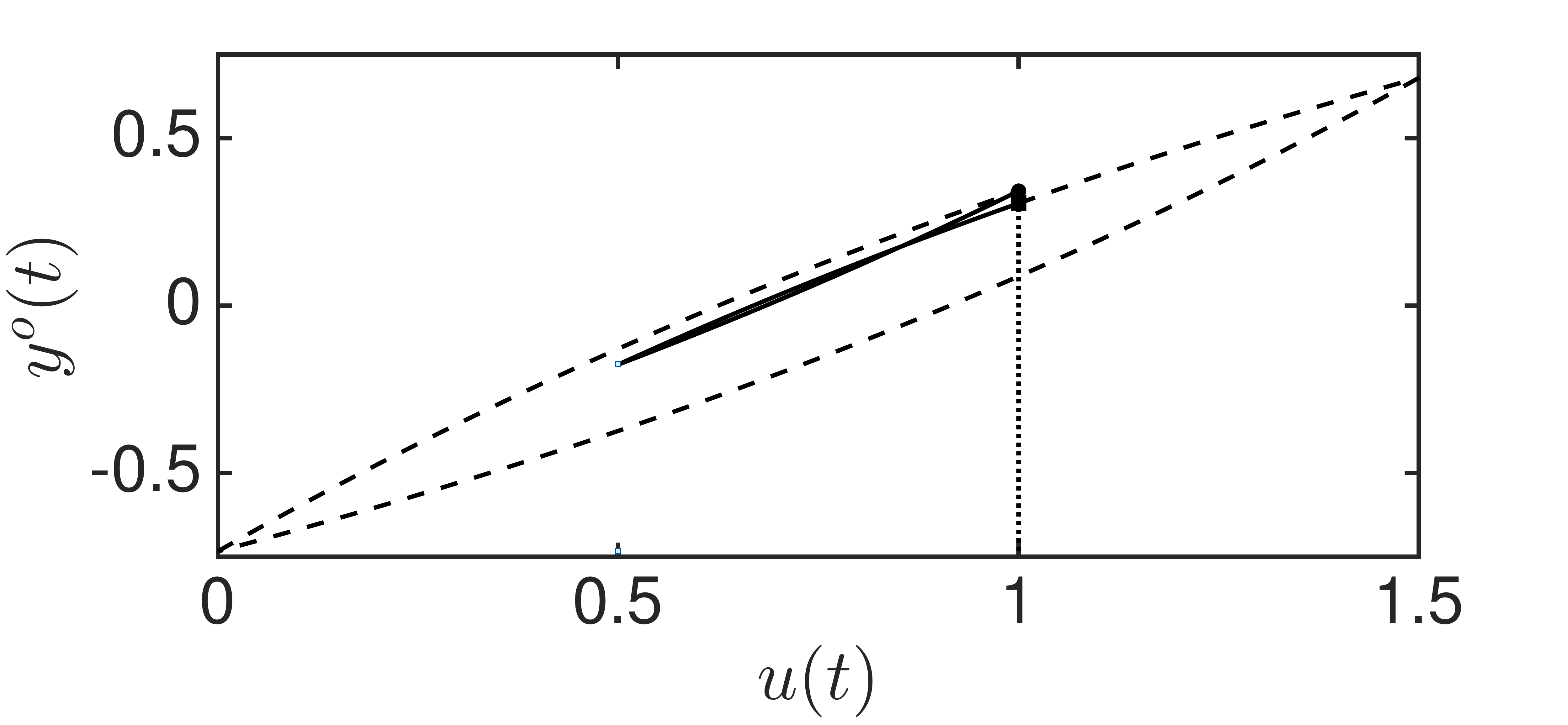}
                \caption{Hysteresis loop $G_u^\circ$ for $\sigma_0=1$ (minor loop in solid line).}
                \label{mlfig8}
\end{figure}

\begin{figure}[H]
                \centering
                \includegraphics[scale=0.25]{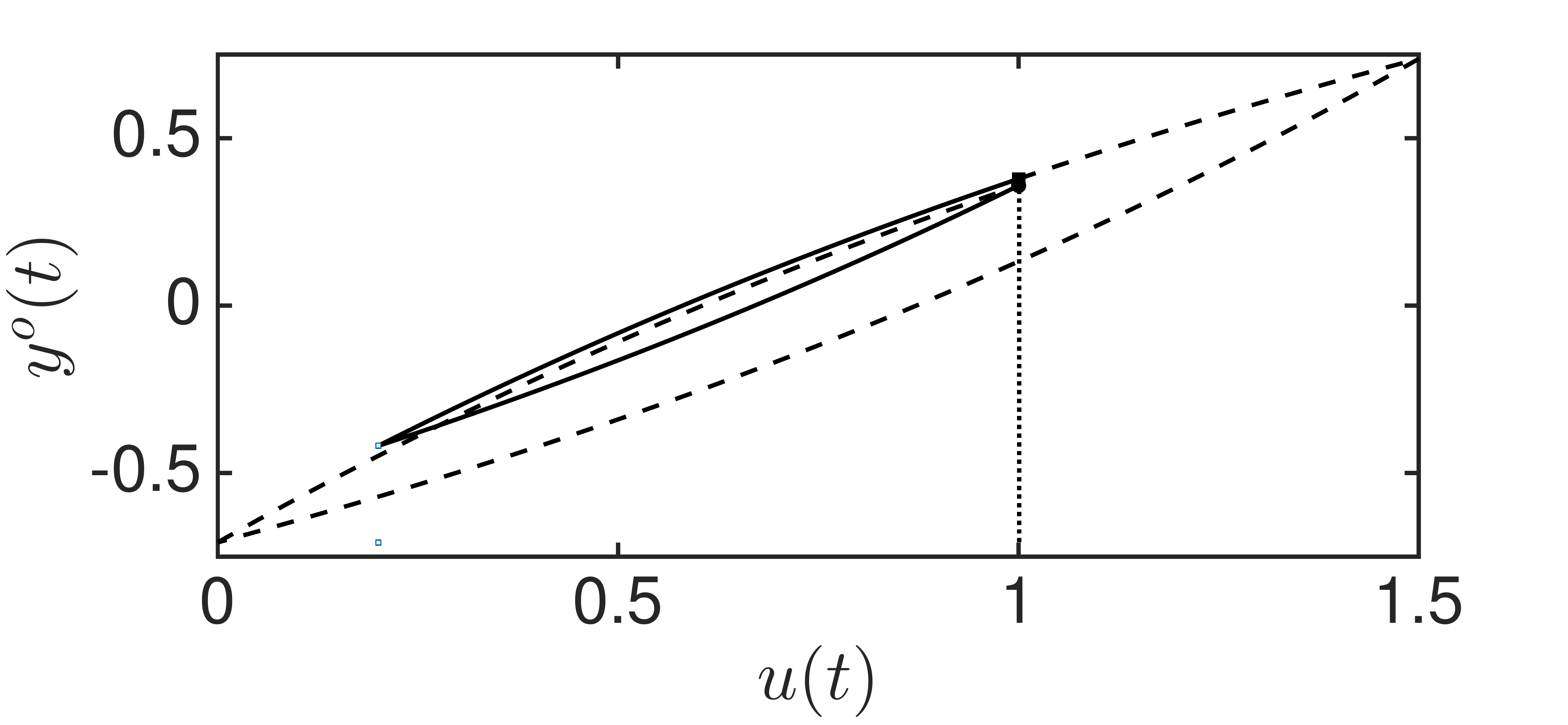}
                \caption{Hysteresis loop $G_u^\circ$ for $\sigma_0=1$ and $u_{\min,2}=0.2$ (minor loop in solid line).}
                \label{mlfig9}
\end{figure}

The hysteresis loop which is given in Figure \ref{mlfig7} is obtained using Equations \eqref{ycircdezeroML}--\eqref{[rho3rho4ycirc]}. Observe that the shape of the minor loop depends greatly on the parameters $\sigma_0$, $F_s$, and on the relative values $u_{\min,2}-u_{\min,1}$, $u_{\max,1}-u_{\min,1}$, and $u_{\max,2}-u_{\min,1}$, see Figures \ref{mlfig8} and  \ref{mlfig9}.


\section{Conclusions} \label{Conclusions and comments}
Although the phenomenon of hysteresis has been studied since the second half of the 19th century, the behavior of minor loops as a specific issue did not emerge as a research subject until the second half of the 20th century. The present paper is framed within the increasing interest in the study of the behavior of minor loops. The originality of this work comes from  being the first to provide an explicit analytic expression of the minor loops of the LuGre and the Dahl models. Our construction can be generalized to  
multimodal input functions giving rise to hysteresis loops with many minor loops. The obtained analytic expressions have been illustrated by means of numerical simulations.

\section*{Acknowledgment}
Funding: The authors are supported by the 
Ministry of Economy, Industry and Competitiveness -- State Research Agency
of the Spanish Government through grant  DPI2016-77407-P (MINECO/AEI/FEDER, UE).
\section*{Compliance with Ethical Standards}
Conflict of Interest:  The authors declare that they have no conflict of interest.


\section*{References}

\begin{thebibliography}{00}

\bibitem{AIRC2012}  N.~Aguirre-Carvajal,  F.~Ikhouane,  J.~Rodellar and  R.~Christenson,  Parametric identification of the Dahl model for large scale MR dampers, Struct. Cont. Health Monit.  19(3)  (2012) 332--347. 

\bibitem{JC08} K.~\AA str\"{o}m, C.~Canudas-de-Wit, Revisiting the LuGre friction model, IEEE Control Syst. Mag. 28(6)  (2008) 101--114.

\bibitem{BM2006}  G.~Bertotti, I.~Mayergoyz, ed.,  The Science of Hysteresis, 3-volume set, Elsevier, Academic Press, Oxford, UK, 2006. 

\bibitem{BS1996}  M.~Brokate, J.~Sprekels, Hysteresis and Phase Transitions, 121, Springer-Verlag, New York, USA, 1996. 

\bibitem{wc95} C.~Canudas de Wit, H.~Olsson, K.~\AA str\"{o}m, P.~Lischinsky,  A new model for control of systems with friction, IEEE Trans.  Autom. Control 40(3) (1995) 419--425.

\bibitem{D1976} P.~Dahl, Solid friction damping of mechanical vibrations, AIAA J. 14(12) (1976) 1675--1682.

\bibitem{FHH2008} R-F.~Fung, C-F.~Han, J-L.~Ha, Dynamic responses of the impact drive mechanism modeled
by the distributed parameter system, Appl. Math. Model. 32(9) (2008) 1734--1743.


\bibitem{GBI2016} I.~Garc\'ia-Ba\~nos, F.~Ikhouane,  A new method for the identification of the parameters of the Dahl model, J. Phys.: Conf. Series.  744(1) (2016) Article ID 012195.

\bibitem{GP2017} C.~Graczykowski, P.~Paw\l owski, Exact physical model of magnetorheological damper, Appl. Math. Model. 47 (2017) 400--424.


\bibitem{HMFA14}  M.~Hamimid, S.~M.~Mimoune, M.~Feliachi, K.~Atallah,  Non centered minor hysteresis loops evaluation based on exponential parameters transforms of the modified inverse Jiles--Atherton model,  Physica B,  451 (2014) 16--19.

\bibitem{I09}  F.~Ikhouane,  Characterization of hysteresis processes, Math. Control Signal Syst.  25(3) (2013) 291--310.

\bibitem{Ikhouane2017} F.~Ikhouane, A survey of the hysteretic Duhem model,  Arch. Comput. Method Eng. 25(4) (2018) 965--1002.


\bibitem{IR2007}  F.~Ikhouane, J.~Rodellar,  Systems with Hysteresis: Analysis, Identification and Control using the Bouc-Wen model, John Wiley \& Sons, The Atrium, Southern Gate, Chichester, England, 2007.

\bibitem{JMS2017} K.~Jankowski, M.~Marsza\l, A.~Stefa\'{n}ski, Formulation of presliding domain non-local memory hysteretic loops based upon modified Maxwell slip model, Tribol. Lett. (2017) 65:56.

\bibitem{Leoni09}  G.~Leoni,  A First Course in Sobolev Spaces,   The American Mathematical Society, USA, 2009.

\bibitem{MNZ93}  J.W.~Macki, P.~Nistri, P.~Zecca,   Mathematical models for hysteresis,  SIAM Rev.  35(1) (1993) 94--123.

\bibitem{M03}  I.~Mayergoyz,  Mathematical Models of Hysteresis, Elsevier Series in Electromagnetism, Elsevier, 2003.


\bibitem{FI13} M.~F.~M.~Naser, F.~Ikhouane, Consistency of the Duhem model with hysteresis, Math. Probl. Eng. (2013) Article ID 586130, 16 pages.


\bibitem{Naser et al.(2015)} M.~F.~M.~Naser, F.~Ikhouane, Hysteresis loop of the LuGre model, Automatica,  59 (2015) 48--53.


\bibitem{OB05}  J.~Oh, D.S.~Bernstein,  Semilinear Duhem model for rate-independent and rate-dependent hysteresis,  IEEE Trans.  Autom. Control 50(5) (2005) 631--645.

\bibitem{Rudin} W.~Rudin, Real and complex analysis, 3rd ed. McGraw-Hill, New York, 1987.


\bibitem{V94}  A.~Visintin, Differential Models of Hysteresis, Springer-Verlag, Berlin, Heidelberg, 1994.

\bibitem{YV1966}  E.H.~Yen, H.R.~Van Der Vaart,   On measurable functions, continuous functions and some related concepts,  Am. Math. Mon. 73(9) (1966) 991--993.


\bibitem{ZRdKAF2017} H.Y.~Zhao, C.~Ragusa, O.~de la Barriere, M.~Khan, C.~Appino, F.~Fiorillo, Magnetic loss versus frequency in non-oriented steel sheets and its prediction: minor loops, PWM, and the limits of the analytical approach, IEEE Trans.  Magn.   53(11) (2017) Article ID 2003804.

\end{thebibliography}
\end{document}